\documentclass[leqno]{amsart}
\usepackage{egastyle}
\usepackage{mymacros}
\usepackage[all,tips]{xy}
\usepackage{url}
\usepackage{bm}
\usepackage{verbatim}

\CompileMatrices

\def\cA{{\mathcal A}}
\def\cB{{\mathcal B}}

\def\cX{{\mathcal X}}

\newcommand{\CC}{\mathbf{C}}

\newcommand{\PP}{\mathbf{P}}

\newcommand{\QQ}{\mathbf{Q}}
\newcommand{\RR}{\mathbf{R}}
\newcommand{\ZZ}{\mathbf{Z}}
\newcommand{\HH}{\mathbf{H}}

\newcommand{\can}{{\operatorname{can}}}

\renewcommand{\setminus}{\smallsetminus}
\newcommand{\MaxSpec}{\operatorname{MaxSpec}}
\newcommand{\la}{{\langle}}
\newcommand{\ra}{{\rangle}}
\newcommand{\an}{{\rm an}}
\newcommand{\s}[1]{{}^{#1}}
\newcommand{\p}{{}'}
\newcommand{\pp}{{}''}
\newcommand{\berk}{{\rm an}}
\newcommand{\Berk}{{\rm an}}
\newcommand{\Gm}{{\mathbf G}_m}

\renewcommand{\Div}{{\rm div}}
\newcommand{\mult}{{\rm mult}}
\newcommand{\rel}{{\rm rel}}
\newcommand{\into}{\hookrightarrow}
\newcommand{\onto}{\twoheadrightarrow}
\def\isomap{{\buildrel \sim\over\longrightarrow}}

\newcommand{\Xhat}{{\hat{X}}}

\newcommand\pu[1]{\{\kern -2pt\{#1\}\kern -2pt\}}

\DeclareMathOperator{\val}{val}

\DeclareMathOperator{\trop}{trop}
\DeclareMathOperator{\Trop}{Trop}
\DeclareMathOperator{\vertices}{vert}
\DeclareMathOperator{\Star}{Star}

\DeclareMathOperator{\inn}{in}

\DeclareMathOperator{\relint}{relint}

\DeclareMathOperator{\dist}{dist}

\newcommand\T{\bT}
\newcommand\A{\bA}

\newcommand\G{\bG}
\renewcommand\top{\mathrm{top}}

\title{Nonarchimedean geometry, tropicalization, and metrics on curves}

\author{Matthew Baker} 
\email{mbaker@math.gatech.edu}
\address{School of Mathematics, Georgia Institute of Technology, Atlanta GA 30332-0160, USA}
\author{Sam Payne}
\email{sam.payne@yale.edu}
\address{Mathematics Department, Yale University, New Haven, CT 06511}
\author{Joseph Rabinoff}
\email{rabinoff@math.harvard.edu}
\address{Department of Mathematics, Harvard University, Cambridge, MA 02138}

\begin{document}

\begin{abstract}
  We develop a number of general techniques for comparing analytifications
  and tropicalizations of algebraic varieties.  Our basic results include
  a projection formula for tropical multiplicities and a generalization of
  the Sturmfels-Tevelev multiplicity formula in tropical elimination
  theory to the case of a nontrivial valuation.  For curves, we explore in
  detail the relationship between skeletal metrics and lattice lengths on
  tropicalizations and show that the maps from the analytification of a
  curve to the tropicalizations of its toric embeddings stabilize to an
  isometry on finite subgraphs.  Other applications include
  generalizations of Speyer's well-spacedness condition and the
  Katz-Markwig-Markwig results on tropical $j$-invariants.
\end{abstract}

\maketitle

\vspace{-10 pt}

\begin{rem*}
  Note that the literature contains a number of references to an earlier
  preprint version of this paper,
  \url{arxiv:1104.0320v2}.
  This final version differs from that one in several respects.  The numbering
  (of equations, paragraphs, sections, theorems, etc.) has changed.  The former
  Section~5, \textit{The structure theory of analytic curves}, was extracted and
  published separately~\cite{bpr:berk_curves}.  Furthermore, much expository
  material and many examples in the remaining sections have been omitted.  The
  earlier preprint version remains available on the arXiv.
\end{rem*}

\section{Introduction}

The recent work of Gubler \cite{gubler:bogomolov,gubler:tropical}, in
addition to earlier work of Bieri-Groves \cite{bieri_groves:valuations},
Berkovich
\cite{berkovich:analytic_geometry,berkovich:locallycontractible1,berkovich:locallycontractible2},
and others, has revealed close connections between nonarchime\-dean
analytic spaces (in Berkovich's sense) and tropical geometry.  One such
connection is given by the second author's theorem that `analytification
is the inverse limit of all tropicalizations' (see
Theorem~\ref{thm:paynehomeo} below).  This result is purely topological,
providing a natural homeomorphism between the nonarchimedean
analytification $X^{\an}$ of a quasiprojective variety $X$ and the inverse
limit of all `extended tropicalizations' of $X$ coming from closed
immersions of $X$ into quasiprojective toric varieties that meet the dense
torus. In this paper, we develop a number of general techniques for
comparing finer properties of analytifications and tropicalizations of
algebraic varieties and apply these techniques to explore in detail the
relationship between the natural metrics on analytifications and
tropicalizations of curves. The proofs of our main results rely on the
geometry of formal models and initial degenerations as well as Berkovich's
theory of nonarchimedean analytic spaces.

\medskip

Let $K$ be an algebraically closed field that is complete with respect to
a nontrivial nonarchime\-dean valuation $\val:K\to\R\cup\{\infty\}$. Let
$X$ be a nonsingular curve defined over $K$.  The underlying topological
space of $X^{\an}$ can be endowed with a `polyhedral' structure locally
modeled on an $\RR$-tree.  The leaves of $X^\an$ are the $K$-points, together
with the `type-4 points' in Berkovich's classification \parref{par:types}.
The non-leaves are exactly those points that are contained in an embedded
open segment, and the space $\HH_\circ(X^\an)$ of non-leaves carries a
canonical metric which, like the polyhedral structure, is defined using
semistable models for $X$.  Our primary reference for these results
is~\cite{bpr:berk_curves}; see also~\cite[\S{4}]{berkovich:analytic_geometry},
\cite{thuillier:thesis}, and~\cite[\S{5}]{baker:aws}.

\medskip

Suppose $X$ is embedded in a toric variety $Y_{\Delta}$ and meets the
dense torus $\T$.  The tropicalization $\Trop(X \cap \T)$ is a
1-dimensional polyhedral complex with no leaves in the real vector space
spanned by the lattice of one parameter subgroups of $\T$.  All edges of
$\Trop(X \cap \T)$ have slopes that are rational with respect to the
lattice of one parameter subgroups, so there is a natural metric on
$\Trop(X \cap \T)$ given locally by lattice length on each edge, and
globally by shortest paths.  The metric space $\HH_\circ(X^\an)$ naturally
surjects onto $\Trop(X \cap \T)$, but this map is far from being an
isometry since infinitely many embedded segments in $\HH_\circ(X^\an)$ are
contracted.  Furthermore, even when an edge of $\HH_\circ(X^\an)$ maps
homeomorphically onto an edge of $\Trop(X \cap \T)$, this homeomorphism
need not be an isometry; see \parref{par:examples} below.  Nevertheless,
each embedded subgraph in $\HH_\circ(X^\an)$ maps isometrically onto its
image in all `sufficiently large' tropicalizations.

\begin{thm} \label{thm:FirstMain} 
  Let $\Gamma$ be a finite embedded
  subgraph in $\HH_\circ(X^\an)$.  Then there is a closed embedding of $X$
  into a quasiprojective toric variety such that $X$ meets the dense torus
  and $\Gamma$ maps isometrically onto its image in $\Trop(X \cap \T)$.
  Furthermore, the set of all such embeddings is stable and hence cofinal
  in the system of all embeddings of $X$ into quasiprojective
  toric varieties whose images meet $\T$.
\end{thm}

\medskip

\noindent Here if $\iota: X\inject Y_\Delta$ and 
$\iota':X\inject Y_{\Delta'}$ are closed embeddings into quasiprojective
toric varieties such that $X$ meets the dense tori $\T$ and $\T'$, then we
say that $\iota'$ \emph{dominates} $\iota$ and we 
write $\iota'\geq\iota$ if there exists an equivariant morphism of toric
varieties $\psi:Y_{\Delta'}\to Y_\Delta$ such that
$\psi\circ\iota'=\iota$ (see~\parref{par:order.embeddings}).  In this case
we have an induced map 
$\Trop(X\cap\T')\to\Trop(X\cap\T)$; the above theorem says in particular
that if $\Gamma$ maps isometrically onto its image
in $\Trop(X\cap\T)$, then the same is true for
$\Gamma\to\Trop(X\cap\T')$.  In other words, the maps from
$\HH_\circ(X^\an)$ to the tropicalizations of toric embeddings of $X$
stabilize to an isometry on every finite subgraph.  

Both the analytification and the tropicalization constructions described
above for subvarieties of tori globalize in natural ways.  The
analytification functor extends to arbitrary finite type $K$-schemes (see
\cite[Chapters~2 and 3]{berkovich:analytic_geometry} or
\cite{berkovich:etalecohomology}), and tropicalization extends to closed
subvarieties of toric varieties as follows.  If $\Delta$ is a fan in
$N_{\RR}$ and $Y_{\Delta}$ is the associated toric variety, then there is
a natural `partial compactification' $N_{\RR}(\Delta)$ of $N_{\RR}$ which
is, set-theoretically, the disjoint union of the tropicalizations of all
torus orbits in $Y_{\Delta}$.  The topology on $N_\R(\Delta)$ is such that
the natural map from $Y_\Delta(K)$ extends to a continuous, proper, and
surjective map $\trop: Y_\Delta^\an \rightarrow N_\R(\Delta)$.  As in the
case where $Y_\Delta$ is the torus $\T$, the tropicalization $\Trop(X)$ of
a closed subvariety $X$ in $Y_\Delta$ is the closure of $\trop(X(K))$ in
$N_\R(\Delta)$, and the extended tropicalization map extends to a
continuous, proper, surjective map from $X^\an$ onto $\Trop(X)$.  See
\cite{payne:analytification, jdr:trop_ps} and \parref{par:extendedtrop}
below for further details.

\begin{thm}[Payne]
\label{thm:paynehomeo}
Let $X$ be an irreducible quasiprojective variety over $K$. Then the
inverse limit of the extended tropicalizations $\Trop(\iota(X))$ over all
closed immersions $\iota : X \into Y_{\Delta}$ into quasiprojective toric
varieties is canonically homeomorphic to the analytification $X^{\an}$.
\end{thm}

\medskip
The inverse limit in Theorem~\ref{thm:paynehomeo} can be restricted to
those closed immersions $\iota$ whose images meet the dense torus
$T_\iota$, and then the homeomorphism maps $X^\an \smallsetminus X(K)$
homeomorphically onto the inverse limit of the ordinary tropicalizations
$\Trop(\iota(X) \cap T_\iota)$.  

\medskip

When $X$ is a curve, our
Theorem~\ref{thm:FirstMain} says that the metric structures on
$\trop(\iota(X) \cap T_\iota)$ stabilize to a metric on the subset $\HH_\circ(X^\an)$ of the inverse limit,
and the restriction of this homeomorphism is an {\em isometry}.
In general, each sufficiently small segment $e$ in $\HH_\circ(X^{\an})$ is
mapped via an affine linear transformation with integer slope onto a
(possibly degenerate) segment $e'$ in $\Trop(X)$.  We write $m_\rel(e)$
for the absolute value of the slope of this map, so if $e$ has length
$\ell$ then its image $e'$ has lattice length $m_\rel(e) \cdot \ell$.  In
Corollary~\ref{cor:Berkovichmultiplicities}, we relate these `expansion
factors' to {\em tropical multiplicities} of edges in $\Trop(X)$.  The
notation is meant to suggest that $m_\rel(e)$ may be thought of in this
context as the relative multiplicity of $e$ over $e'$.  By definition, the
tropical multiplicity $m_{\Trop}(e')$ of an edge $e'$ in a suitable
polyhedral structure on $\Trop(X)$ is the number of irreducible
components, counted with multiplicities, in the {\em initial degeneration}
$\inn_w(X \cap \T)$ for any $w$ in the relative interior of $e'$.  These
tropical multiplicities are fundamental invariants in tropical geometry
and play a key role in the {\em balancing formula}. 
See~\parref{par:tropicalization} for a definition of the
initial degeneration $\inn_w(X \cap \T)$ and further discussion of
tropical multiplicities.

\begin{thm}  \label{thm:Subdivision}
There is a polyhedral structure on $\Trop(X \cap \T)$ with the following properties.
\begin{enumerate}
\item For each edge $e'$ in $\Trop(X \cap \T)$, there are finitely many
  embedded segments $e_1, \ldots, e_r$ in $\HH_\circ(X^\an)$ mapping
  homeomorphically onto $e'$.

\item Any embedded segment in the preimage of $e'$ that is disjoint from
  $e_1 \cup \cdots \cup e_r$ is contracted to a point.

\item The tropical multiplicity of $e'$ is the sum of the corresponding
  expansion factors
  \[
  m_{\Trop}(e') = m_\rel(e_1) + \cdots + m_\rel(e_r).
  \]
\end{enumerate}
\end{thm}

\medskip

\noindent The properties above are preserved by subdivision, so they hold
for any sufficiently fine polyhedral structure on $\Trop(X \cap \T)$.
See Proposition~\ref{prop:compatible.polyhedral} and
Corollary~\ref{cor:Berkovichmultiplicities}. 

\medskip

The tropical multiplicity formula in the above theorem gives an important
connection to nonarchimedean analytic spaces that is not visible from the
definitions.  The formula shows, for example, that if $e'$ is a small
segment in $\Trop(X \cap \T)$ whose tropical multiplicity is equal to $1$,
then there is a unique segment $e$ in $\HH_\circ(X^{\an})$ mapping
homeomorphically onto $e'$, and the length of $e$ is equal to the tropical
length of $e'$.  It is well known that the skeleton of the analytification
of an elliptic curve with bad reduction is a loop of length equal to minus
the valuation of the $j$-invariant (see for example \cite[Remark 4.24]{bpr:berk_curves}), so
these formulas explain earlier results of Katz, Markwig and Markwig on
tropical $j$-invariants of genus one curves in toric surfaces
\cite{KatzMarkwigMarkwig08, KatzMarkwigMarkwig09}.  See, for instance,
Example~\ref{eg:Example2a}.  The following theorem also provides natural
generalizations for genus one curves in higher dimensional toric
varieties, as well as curves of arbitrary genus.

\begin{thm}  \label{thm:Mult1Subgraph}
  Let $\Gamma'$ be a finite embedded subgraph of $\Trop(X \cap \T)$ and
  suppose $\inn_w(X \cap \T)$ is irreducible and generically reduced for
  every $w$ in $\Gamma'$.  Then there is a unique embedded subgraph
  $\Gamma$ in $\HH_\circ(X^\an)$ mapping homeomorphically onto $\Gamma'$,
  and this homeomorphism is an isometry.
\end{thm}

\medskip

\noindent
See \S\ref{section:ellipticcurves} for details on deducing the tropical
$j$-invariant results of Katz, Markwig and Markwig from the above theorem.

\medskip

The expansion factors $m_{\rel}(e)$ in our tropical multiplicity formula
are often computable in practice.  If $X$ is an affine curve embedded in
the torus $\G_m^n$ via an $n$-tuple of invertible regular functions
$f_1,\ldots,f_n$, then
\[ m_{\rel}(e) = \gcd(s_1(e),\ldots,s_n(e)), \] where $s_i(e)$ is the
absolute value of the slope of the integer-affine function $\log |f_i|$
along the edge $e$.  See Remark~\ref{rem:mult.is.gcd}.  The quantities
$s_i(e)$ are easily calculated from the divisors of $f_1,\ldots,f_n$ using
the `Slope Formula' of~\cite[Theorem~5.15]{bpr:berk_curves}.

\medskip

In concrete situations, it is useful to be able to certify that a given
tropicalization map faithfully represents a large piece of the nonarchimedean analytification $X^{\an}$ 
(e.g. the `minimal skeleton' $\Sigma$ of $X^{\an}$ in the
sense of Berkovich \cite{berkovich:analytic_geometry} or
\cite[Corollary~4.23]{bpr:berk_curves}) using only 
`tropical' computations (e.g. Gr{\"o}bner complex computations which have
been implemented in computer algebra packages such as Singular or
Macaulay2), as opposed to calculations with formal models that have not
been implemented in a systematic way in any existing software package.  We
prove that a tropicalization map represents $\Sigma$ faithfully, meaning
that the map is an isometry on $\Sigma$, provided that certain
combinatorial and topological conditions are satisfied.  
Our results on faithful representations are presented in conjunction with some observations about
initial degenerations which help explain the special role played by trivalent
graphs in the literature on tropical curves (cf. Theorem~\ref{thm:combmult1theorem} and
Remark~\ref{rem:combmult1remark}).

\medskip

We explore tropicalizations of elliptic curves in detail as a concrete
illustration of our methods and results.  We are able to say some rather
precise things in this case; for example, we show that every elliptic
curve $E/K$ with multiplicative reduction admits a closed embedding in
$\PP^2$ whose tropicalization faithfully (and certifiably) represents the
minimal skeleton of $E^{\an}$.  
Furthermore, we interpret Speyer's `well-spacedness
condition' for trivalent tropicalizations of totally degenerate genus 1
curves \cite{speyer:uniformizing} as a statement about rational functions
on the analytification of the curve, and prove generalizations of this
condition for nontrivalent tropicalizations, and for genus 1 curves with
good reduction.

\medskip

The paper concludes with a generalization of the important Sturmfels-Tevelev multiplicity formula,
which calculates $\Trop(\alpha(X))$ (as a weighted polyhedral complex) in terms
of $\Trop(X)$ when $\alpha : \T \to \T'$ is a homomorphism of tori which
induces a generically finite map from a subvariety $X$ in $\T$ onto its
image.  The 
multiplicity formula in \cite{sturmfels_tevelev:elimination} is formulated
and proved in the `constant coefficient' setting, where $K$ is the field
of Puiseux series over an algebraically closed coefficient field $k$ of
characteristic $0$ and $X$ is defined over $k$.  We use the methods of this paper to generalize the
Sturmfels-Tevelev formula to the case where $X$ is any closed subvariety
of a torus $\T$ defined over a complete and algebraically closed
nonarchimedean field $K$.
 
\medskip

Philosophically speaking, there are at least two long-term goals to this
paper.  On the one hand, we believe that the systematic use of modern
tools from nonarchimedean geometry is extremely useful for understanding
and proving theorems in tropical geometry.  This paper takes several steps
in that direction, establishing some new results in tropical geometry via
Berkovich's theory and the Bosch-L{\"u}tkebohmert-Raynaud theory of
admissible formal schemes.  On the other hand, much of this paper can be
viewed as a comparison between two different ways of approximating
nonarchimedean analytic spaces.  Nonarchimedean analytic spaces have
proved to be useful in many different contexts, but the topological spaces
underlying them are wildly branching infinite complexes which are
difficult to study directly, so one usually approximates them with finite
polyhedral complexes.  One such approximation goes through skeleta of nice
(e.g. semistable) formal models (cf.\ \cite[Theorem 5.2]{bpr:berk_curves}),
another through (extended)
tropicalizations (cf. Theorem~\ref{thm:paynehomeo}).  Our
Theorem~\ref{thm:MainThm4} shows that, in the case of curves, these two
approximations have the same metric structure in the limit (though the
metrics may be different at any given finite level).

\medskip

For further details and examples, we refer the reader to the expanded version of
this paper at \texttt{arXiv:1104.0320}.  Since this paper was written, there
have been a large number of follow-up articles: the reader may also want to
read~\cite{GRW}, in which many of the results in this article are extended to
higher dimensions; \cite{CHW,dp,chan_sturmfels:honeycomb}, in which several
interesting examples of faithful tropicalizations are given; and
\cite{abbr14:lifting_harmonic_morphism_I,abbr14:lifting_harmonic_morphism_II},
in which ``relative'' versions of some of the results in this paper are used to
prove tropical lifting theorems.

\medskip

\noindent \textbf{Acknowledgments.}  \thanks{The authors would like to
  express their thanks to Mar{\'i}a Ang{\'e}lica Cueto, Eric Katz, Brian
  Osserman, David Speyer, and Josephine Yu for helpful discussions, and to
  Bernd Sturmfels for his interest and encouragement.  Special thanks are
  due to Walter Gubler for some illuminating conversations which helped
  shape the direction of this work.  The authors also thank Melanie Dunn
  for computing several nice examples of tropicalizations, and the
  anonymous referees for many helpful comments.  The first
  and second authors were supported in part by NSF Research Grants DMS-0901487 and DMS-1068689, respectively.}

\section{Basic notions and examples}

Here we give a brief overview of the basic notions necessary to understand
the theorems stated in the introduction, followed by a few key examples
illustrating these results.  Throughout this paper, $K$ is an
algebraically closed field that is complete with respect to a nontrivial
nonarchimedean valuation
\[
\val:K\to\R\cup\{\infty\}.
\]
We let $G = \val(K^\times)$ be its value group, $R = \val\inv([0,\infty])$
its valuation ring, $\fm\subset R$ the maximal ideal, and $k = R/\fm$ its
residue field (which is algebraically closed by 
\cite[\S{2.1} Proposition 3]{robert:padic_analysis}).
Let $|\cdot| = \exp(-\val(\cdot))$ be the absolute value on $K$ associated
to the valuation.

\medskip

\paragraph[Tropicalization]
\label{par:tropicalization}
Let $M$ be a free abelian group of rank $n$, let $\T = \Spec(K[M])$ be the
$K$-torus with character group $M$, and let $N = \Hom(M,\ZZ)$ be the dual
lattice.  If $X$ is a closed subscheme of $\T$, there is a natural 
{\em tropicalization map}
\[
\trop : X(K) \to N_{\RR},
\]
where $N_{\RR} = \Hom(M,\R)$.  The image of a point $x$ in $X(K)$ is the
linear function taking $u \in M$ to the valuation of the corresponding
character evaluated at $x$.  Then $\Trop(X)$ is the closure of
$\trop(X(K))$ in the Euclidean topology on $N_{\RR}$.  Note that the
choice of an isomorphism $M \cong \ZZ^n$ induces an identification of $\T$
with $\G_m^n$.  In such coordinates, the tropicalization map sends a point
$(x_1,\ldots,x_n)$ in $X(K)$ to $(\val(x_1),\ldots,\val(x_n))$ in $\RR^n$.

One of the basic results in tropical geometry says that if $X$ is an
integral subscheme of $\T$ of dimension $d$ then $\Trop(X)$ is the
underlying set of a connected `balanced weighted integral $G$-affine
polyhedral complex' of pure dimension $d$.  We do not define all of these
terms here, but briefly recall how one gets a polyhedral complex and
defines weights on the maximal faces of this complex.  Let $w$ be a point in $N_G = \Hom(M,G)$.  
The `tilted group ring' $R[M]^w$ is the subring of $K[M]$ consisting of 
Laurent polynomials $a_1 x^{u_1} + \cdots + a_r x^{u_r}$ such that
\[  
\val(a_i) + \la u_i,w \ra \geq 0
\]
for all $i$.
The $R$-scheme $\T^w = \Spec R[M]^w$ is
a torsor for the torus $\Spec R[M]$, and its generic fiber is canonically
isomorphic to $\T$.  If $X$ is a closed subscheme of $\T$ defined by an
ideal $\fa\subset K[M]$ then
\[
X^w = \Spec \big(R[M]^w/\left( \fa \cap R[M]^w \right) \big)
\]
is a flat $R$-scheme with generic fiber $X$, which we call the 
{\em tropical integral model}
associated to $w$.  It is exactly the closure of $X$ in $\T^w$.  The
special fiber $\inn_w(X)$ of $X^w$ is called the {\em initial
  degeneration} of $X$ with respect to $w$ and is the subscheme of the
special fiber of $\T^w$ cut out by the $w$-initial forms of Laurent
polynomials in $\fa$, in the sense of generalized Gr{\"o}bner theory.

The scheme $\T^w$ is not proper, so points in $X(K)$ may fail to have
limits in the special fiber.  Indeed, the special fiber $\inn_w(X)$ is
often empty.  One of the fundamental theorems in tropical geometry says
that $w$ is in $\Trop(X)$ if and only if $\inn_w(X)$ is not empty.%
\footnote{In the special case where $\T$ has dimension one and $X$ is the
  zero locus of a Laurent polynomial $f$, this is equivalent to the
  statement that $f$ has a root with valuation $s$ if and only if $-s$ is
  a slope of the Newton polygon of $f$.}  Moreover, $\Trop(X)$ can be
given the structure of a finite polyhedral complex in such a way that
whenever $w$ and $w'$ belong to the relative interior of the same face,
the corresponding initial degenerations $\inn_w(X)$ and $\inn_{w'}(X)$ are
$\T$-affinely equivalent.

We define the {\em multiplicity} $m_{\Trop}(w)$ of a point $w$ in
$\Trop(X)$ to be the number of irreducible components of $\inn_w(X)$,
counted with multiplicities.  In particular $m_{\Trop}(w) = 1$ if and only
if $\inn_w(X)$ is irreducible and generically reduced.  These tropical
multiplicities are constant on the relative interior of each face $F$ of
$\Trop(X)$, and we define the multiplicity $m_{\Trop}(F)$ to be
$m_{\Trop}(w)$ for any $w$ in the relative interior of $F$.  The
multiplicities for maximal faces are the `weights' mentioned above that
appear in the balancing condition.  These weights have the following
simple interpretation for hypersurfaces.

\begin{rem} \label{rem:newton.complex}
  If $X = V(f)$ is a hypersurface then $\Trop(X)$ is the corner locus of the
  convex piecewise linear function associated to a defining
  equation $f$~\cite[Section 2.1]{ekl:non_arch_amoebas}.  In this case,
  $\Trop(X)$ has a unique minimal polyhedral structure, and the initial
  degenerations are essentially constant on the relative interior of each
  face.  There is a natural inclusion reversing bijection between the
  faces of $\Trop(X)$ in this minimal polyhedral structure and the
  positive dimensional faces of the \emph{Newton polytopal complex} (or
  \emph{Newton complex}) of $f$: a face of
  $\Trop(X)$ corresponds to the convex hull of the monomials whose
  associated affine linear function is minimal on that face.  In
  particular, the maximal faces of $\Trop(X)$ correspond to the edges of
  this Newton complex.  In this special case, the multiplicity of a
  maximal face is the lattice length of the corresponding
  edge.\footnote{This is explained in~\cite[Example
    3.16]{sturmfels_tevelev:elimination} in the special case where $X$ is
    irreducible, $K$ is the field of Puiseux series over $k$, and $X$ is
    defined over $k$.  The arguments given there work in full generality.}
  The relationship between the tropical hypersurface and the Newton
  complex is also explained in more detail in~\cite[\S8]{jdr:trop_ps}.
\end{rem}

\paragraph[Analytification] \label{par:analytification} Let $A$ be a
finite-type $K$-algebra.  The \emph{Berkovich spectrum} of $A$, denoted
$\sM(A)$, is defined to be the set of multiplicative seminorms $\|\cdot\|$
on $A$ extending the absolute value $|\cdot|$ on $K$.  
The Berkovich spectrum $\sM(A)$ is the underlying set of the
nonarchimedean analytification $X^\an$ of $X = \Spec(A)$.  The topology on
$X^\an$ is the coarsest such that the map $\|\cdot\|\mapsto \|f\|$ is
continuous for every $f\in A$; this coincides with the subspace topology
induced by the inclusion of $X^\an$ in $\R^A$.

\begin{rem*}
  We will often write $\A^1_\an$ for $\A^{1,\an}$ and $\bP^1_\an$ for 
  $\bP^{1,\an}$, etc. 
\end{rem*}

If $X$ is connected then $X^\an$ is a path-connected locally compact
Hausdorff space that naturally contains $X(K)$ as a dense subset; a point
$x \in X(K)$ corresponds to the seminorm $\|\cdot\|_x$ given by $\|f\|_x =
|f(x)|$.  The analytification procedure $X \mapsto X^{\an}$ gives a
covariant functor from the category of locally finite-type $K$-schemes to
the category of topological spaces.\footnote{The analytification $X^\an$
  has additional structure, including a structure sheaf, and $X \mapsto
  X^\an$ may also be seen as a functor from locally finite-type
  $K$-schemes to locally ringed spaces.  See
  \cite[\S{2.3},\S{3.1},\S{3.4}]{berkovich:analytic_geometry} for more details.}

If $X$ is a closed subvariety of $\T$ then the tropicalization map
described in the previous section extends from $X(K)$ to a continuous and
proper map
\[
\trop: X^\an \rightarrow N_\R
\]
taking a seminorm $\| \cdot \|$ to the linear function 
$u \mapsto -\log \|x^u \|$, and the image of this map is exactly
$\Trop(X)$.  In other words, 
$\trop(X^\an)$ is the closure of the image of $X(K)$ in $N_\R$.

\paragraph[Metric structure of analytic curves]
\label{par:metricstructure}
There is a natural metric on $\PP^1_\an \smallsetminus \PP^1(K)$, see \cite{bpr:berk_curves} for details.
We write $\HH(\PP^1_\an)$ to denote
$\PP^1_\an \smallsetminus \PP^1(K)$ with this metric
structure.\footnote{It is important to note that the metric topology on
  $\HH(\PP^1_\an)$ is much finer than the subspace topology on $\PP^1_\an
  \smallsetminus \PP^1(K)$.  Our notation follows \cite{baker_rumely:book}
  and reflects the fact that the metric on $\HH(\PP^1_\an)$ is
  $0$-hyperbolic in the sense of Gromov.}  

The metric on $\HH(\PP^1_\an)$ has the important property that, roughly
speaking, $\log |f|$ is piecewise affine with integer slopes for any
nonzero rational function $f \in K(T)$.  More precisely, suppose that $f$ is
nonconstant and that $\widehat \Sigma$ is the minimal closed connected subset
of $\PP^1_\an$ containing the set $S$ of zeros and poles of $f$.  Let
$\Sigma= \widehat \Sigma \smallsetminus S$.  Then:
\begin{enumerate}
\item The subspace $\Sigma$ of $\HH(\PP^1_\an)$ is a metric graph with finitely many edges, in which the edges whose closures meet $K$ have infinite length.
\item The restriction of $\log|f|$ to $\Sigma$ is piecewise affine with integer slopes. 
\item There is a natural retraction map from $\PP^1_\an$ onto $\widehat \Sigma$.
\item The function $\log|f|$ from $\PP^1_\an$ to $\RR \cup \{ \pm \infty\}$ factors through the retraction onto 
$\widehat \Sigma$, and hence is determined by its restriction to $\Sigma$.
\end{enumerate}
The metric on the complement of the set of $K$-points in the analytification of an arbitrary algebraic curve is induced by the metric on $\PP^1_\an \smallsetminus \PP^1(K)$ via semistable decomposition; see \cite{bpr:berk_curves} for details.

There is also a notion of a skeleton of a smooth and connected but not
necessarily complete curve $X$.  Let $\hat X$ be the smooth
compactification of $X$ and let $D = \hat X\setminus X$ be the set of
`punctures'.  Choose a semistable model $\cX$ of $\hat X$ such that the
punctures reduce to distinct smooth points of the special fiber
$\bar\cX$.  Then there is unique minimal closed connected subset $\Sigma$
of $X^\an$ containing the skeleton $\Sigma_\cX$ of $\hat X$ and whose
closure in $\hat X^\an$ contains $D$.  We call $\Sigma$ the \emph{skeleton} of $X$
associated to $\cX$.  As above, there is a canonical retraction map
$\tau_\Sigma: X^\an\surject\Sigma$.  If $X\subset\T$ then the
tropicalization map $\trop:X\to N_\R$ factors through $\tau_\Sigma$.
There is a skeleton which is minimal over all models
$\cX$ if $2 - 2g(\hat X) - \#D\leq 0$.  See \cite{bpr:berk_curves} for a complete discussion of the skeleta of a curve.

\paragraph[Examples]
\label{par:examples}
To illustrate our main results concerning the relationship between
analytification and tropicalization in the case of curves, taking into
account the metric structure on both sides, we present the following
examples.  In each example, we fix a specific coefficient field for concreteness.

Our 
first example shows how a loop in the analytification of a genus 1
curve can be collapsed onto a segment of multiplicity greater than 1.

\begin{eg}
\label{eg:Example2b}
Let $K$ be the completion of the field of Puiseux series $\C\pu t$.
Consider the genus 1 curve $\hat E\subset\PP^2$ over $K$ defined by the
Weierstrass equation $y^2 = x^3 + x^2 + t^4$, and let $E = \hat E\cap\G_m^2$.
The $j$-invariant of $\hat E$ has valuation $-4$, so $\hat E$ has multiplicative
reduction and the minimal skeleton $\Sigma$ of $\hat E$ is isometric to a
circle of circumference $4$.
In this example, $\Trop(E)$ does not have a cycle even though $\hat E^{\an}$
does; it is interesting to examine exactly what tropicalization is doing
to $\hat E^{\an}$.   Let $\Gamma$ be the minimal skeleton of
$E$; as above, $\trop$ factors through the retraction of $E^{\an}$
onto $\Gamma$. 
Figure~\ref{fig:tropE} shows the restriction of $\trop$ to $\Gamma$.

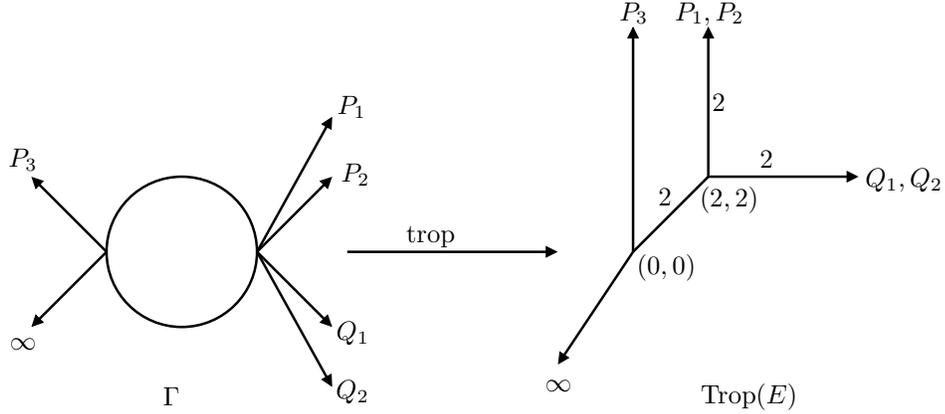
\genericfig[ht]{tropE}{The minimal skeleton $\Gamma\subset E^\an$ and the
  tropicalization $\Trop(E)$ from Example~\ref{eg:Example2b}.  The points
  $P_i,Q_j$ are defined as 
  follows: the rational function $x$ on $\hat E$ has
  divisor $(Q_1)+(Q_2)-2(\infty)$ where $Q_1 = (0,t^2)$ and
  $Q_2=(0,-t^2)$, and $y$ has divisor $(P_1)+(P_2)+(P_3)-3(\infty)$, where 
  $\val(x(P_1)) = \val(x(P_2)) = 2$ and $\val(x(P_3)) = 0$.}

The tropicalization map sends $\Sigma$ $2$-to-$1$ onto its image in
$\Trop(E)$, which is a segment of tropical length 2 and tropical
multiplicity 2.  Locally on $\Sigma$ the tropicalization map is an
isometry.  Each of the rays of $\Gamma$ emanating from $\Sigma$ maps
isometrically onto its image.  The two rays in $\Trop(E)$ with
multiplicity 1 have unique preimages in $\Gamma$, while there are two
distinct rays in $\Gamma$ mapping onto each of the two rays in $\Trop(E)$
of multiplicity 2.  
\end{eg}

\bigskip

The following example illustrates a different kind of collapse, where a segment $e$ in the minimal skeleton of the analytification is collapsed to a point, i.e. the relative multiplicity $m_\rel(e)$ is zero.

\begin{eg}
\label{eg:Example3}
Let $p \geq 5$ be a prime, let $K = \CC_p$, and let $k\cong\bar\F_p$ be its residue
field.  Let $X \subset \Gm^2$ be the affine curve over $K$ defined by the
equation $f(x,y) = x^3y - x^2y^2 - 2xy^3-3x^2y+2xy-p = 0$.  The curve
$\hat{X} \subset \PP^2$ defined by the homogenization $\hat{f}(x,y,z) =
x^3y - x^2y^2 - 2xy^3-3x^2yz+2xyz^2-pz^4 = 0$ is a smooth plane quartic of
genus $3$, and the given equation $\hat{f}(x,y,z)=0$ defines the minimal
regular proper semistable model $\cX$ for $\hat{X}$ over $\QQ_p$.  The
special fiber $\bar{\cX}$ of $\cX$ consists of four (reduced) lines in
general position in $\PP^2_k$, since $\hat{f}$ mod $p$ factors as
$xy(x+y-z)(x-2y-2z)$.  The tropicalization $\Trop(X) \subset \RR^2$
consists of a triangle with vertices $(0,0),(1,0),(0,1)$ together with
three rays emanating from these three vertices in the directions of
$(-1,-1),(3,-1),(-1,3)$ respectively.  The three bounded edges/rays
incident to $(0,0)$ all have tropical multiplicity $2$, and all other bounded
edges/rays in $\Trop(X)$ have tropical multiplicity $1$.  
Let $\Sigma$ be the skeleton $\Sigma_\cX$ of $\hat X$ and let $\Gamma$ be
the minimal skeleton of $X$.  Then $\Sigma$ is a tetrahedron (with
vertices corresponding to the four irreducible components of $\bar{\cX}$)
with six edges of length $1$, because this is the dual graph of a regular semistable model defined over $\ZZ_p$, and $\Gamma$ is obtained from $\Sigma$ by
adding a ray emanating from each vertex of $\Sigma$ toward the zeros and
poles of $x$ and $y$, namely toward the points
$(0:1:0),(1:0:0),(2:1:0),(-1:1:0)$. 
The tropicalization
map $\trop : X^{\an} \to \Trop(X) \subset \RR^2$ factors though the
retraction map $X^{\an} \to \Gamma$.
See Figure~\ref{fig:genus3}.

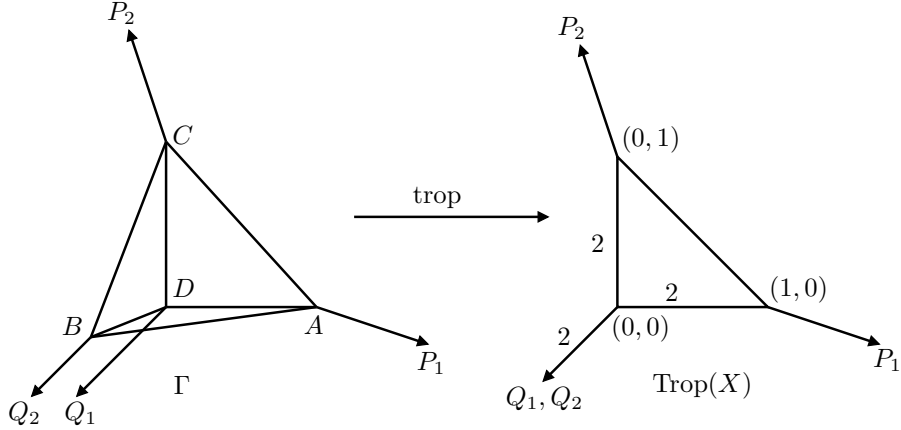
\genericfig[ht]{genus3}{The skeleton $\Gamma\subset X^{\an}$ and the
  tropicalization $\Trop(X)$, where $X$ is the curve from
  Example~\ref{eg:Example3}.  Here  $x$ has divisor $3(P_1)-(P_2)-(Q_1)-(Q_2)$ and $y$ has divisor
  $3(P_2)-(P_1)-(Q_1)-(Q_2)$ on $\Xhat$, where $P_1 = (0:1:0)$, $P_2 = (1:0:0)$,
  $Q_1 = (2:1:0)$, and $Q_2 = (-1:1:0)$.  The points $A,B,C,D\in\Gamma$
  correspond to the irreducible components $x=0, x+y=z, y=0,$ and $x-2y=2z$, 
  respectively, of $\bar{\cX}$.  The collapsed segment is $\overline{BD}$.}

The points of $\Sigma$ corresponding to the irreducible components $x+y=z$ and $x-2y=2z$ of $\bar{\cX}$, as
well as the entire edge of $\Sigma$ connecting these two points, get mapped by $\trop$ to the point $(0,0) \in
\Trop(X)$.  This edge therefore has expansion factor zero with respect to $\trop$.  The other five bounded edges of 
$\Sigma \subset \Gamma$ map isometrically (i.e., with expansion factor one) onto their images in $\Trop(X)$.
In fact, the tropicalization map is a local isometry everywhere on $\Gamma$ except along the bounded edge
which is contracted to the origin.

\end{eg}

\bigskip

Our final example, which is meant to illustrate Theorem~\ref{thm:Mult1Subgraph}, is a genus 1 curve with multiplicative reduction for which the tropicalization map takes the minimal skeleton isometrically onto its image.  This example is also discussed in \cite[Example~5.2]{KatzMarkwigMarkwig09}).

\begin{eg}
\label{eg:Example2a} 
Let $K$ be the completion of the field $\CC\pu t$ of Puiseux series.  
Consider the curve $E'$ in $\Gm^2$ cut out by the equation
$f(x,y) = x^2 y + x y^2 + \frac{1}{t} xy + x + y$.   Its closure in $\PP^2$ is the smooth projective genus 1 curve $\hat E'$ defined by $\hat{f}(x,y,z) = x^2 y + x y^2 +
\frac{1}{t} xyz + xz^2 + yz^2$.  

Using the description of $\Trop(E')$ as the corner locus of the convex
piecewise-linear function associated to $f$, one sees that $\Trop(E')$
consists of a square with side length 2 plus one ray emanating from each
corner of the square; see Figure~\ref{fig:tropEprime}.  By restricting $f$
to faces of the Newton complex (see Remark~\ref{rem:newton.complex}), one checks
that $\inn_w(E')$ is reduced and irreducible for 
every $w$ in $\Trop(E')$.  Therefore, by Theorem~\ref{thm:Mult1Subgraph}
there is a unique graph $\Gamma$ in the analytification of $E'$ mapping
isometrically onto $\Trop(X)$.

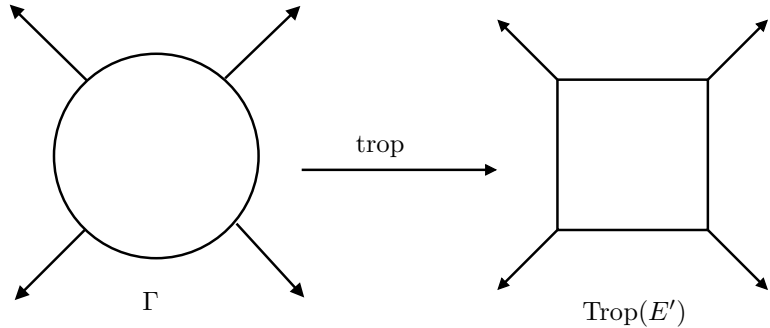
\genericfig[ht]{tropEprime}{The tropicalization of the elliptic curve
  $\hat E'$ from Example~\ref{eg:Example2a}.  The edges in the square each
  have lattice length $2$.} 

In particular, the analytification of $E'$ contains a loop of length 8.  One can check by an explicit computation that $\val(j(\hat E')) = -8$, which is consistent with the
fact that the analytification of a smooth projective genus 1 curve is either contractible (if the curve has good reduction) or else contains a unique loop of length $-\val(j)$
(if the curve has multiplicative reduction).  
\end{eg}

\section{Admissible algebras and nonarchimedean analytic spaces}
\label{section:admissible.R.algebras}

Recall that $K$ is an algebraically closed field 
that is complete with respect to a nontrivial  
nonarchimedean valuation $\val:K\to\R\cup\{\infty\}$.  
Let $|\cdot| = \exp(-\val(\cdot))$ be the associated absolute value.  Let $R$
be the valuation ring of $K$, let $\fm\subset R$ be its maximal ideal, and
let $k = R/\fm$ be its residue field.  Choose a nonzero element
$\varpi\in\fm$ (as $R$ is not noetherian, it has no uniformizer), so $R$
is $\varpi$-adically complete.   
Let $G = \val(K^\times)\subset\R$ be the value group, which is divisible.

In this section we define some notations and collect some results about
admissible $R$-algebras and Raynaud's generic fiber functor that will be needed
in the sequel.  We refer the reader to~\cite{bl:fmI} and
\cite{bosch:lecture_notes} for a detailed discussion of admissible formal
schemes and the Raynaud generic fiber functor.

\paragraph[Admissible formal schemes]
An \emph{admissible $R$-algebra} is a topological $R$-algebra $A$ which is flat
and topologically of finite presentation. A formal scheme $\fX$ over $\Spf(R)$
which is locally isomorphic to the formal spectrum of an admissible $R$-algebra
(in its $\varpi$-adic topology) is called an \emph{admissible formal scheme}.
All admissible formal schemes appearing in this paper will be assumed to be
quasi-compact and separated.  If $\fX = \Spf(A)$ then we write $\fX^\an$ for the
$K$-affinoid Berkovich analytic space $\sM(A\tensor_R K)$.  This construction
globalizes to give a functor from admissible formal schemes to Berkovich
analytic spaces that is compatible with fiber products.  The image of a formal
scheme under this functor is called the \emph{Raynaud generic fiber}.

An admissible formal scheme $\fX$ with reduced special
fiber is called a \emph{formal analytic variety}.  If $\Spf(A)$ is a
formal affine open subset of a formal analytic variety then $A$ is reduced
and $A$ is equal to the full ring of power-bounded elements 
in $A_K \coloneq A\tensor_R K$ by Proposition~\ref{prop:canonical.reduction}.
The canonical reduction $\td A_K$ of $A_K$ therefore coincides with
$A\tensor_R k$, so if $\Spf(B)$ is a formal affine open subset of
$\Spf(A)$ then $\Spec(\td B_K)$ is an affine open subset of 
$\Spec(\td A_K)$.  (Our definition of a formal analytic variety differs
from the original one given in~\cite{bosch_lutk:uniformization}, but the above
argument shows that the two definitions are equivalent in our situation.) 

\begin{notn*}
  Let $\fX$ be an $R$-scheme or a formal $R$-scheme.  We denote the
  special fiber $\fX\tensor_R k$ of $\fX$ by $\bar\fX$.
\end{notn*}

\paragraph[Reductions of analytic spaces]
Let $\cA$ be a $K$-affinoid algebra.  We denote by
$\mathring\cA$, $\check\cA$, and $\td\cA = \mathring\cA/\check\cA$ the subring
of power-bounded elements,
the ideal of topologically nilpotent elements, and the canonical reduction,
respectively.  Setting $\sX=\sM(\cA)$, there is a
\emph{reduction map} $\red : \sX \to \td\sX := \Spec(\td\cA)$.
By~\cite[Corollary~2.4.2 and Proposition~2.4.4]{berkovich:analytic_geometry}
this map is surjective and
anti-continuous, in the sense that the inverse image of an open subset is
closed.  Similarly, if $\fX$ is an admissible formal scheme over $\Spf(R)$ then 
there is a canonical surjective and anti-continuous \emph{reduction map}
$\red: \fX^\an\to\bar\fX$ which
coincides with the map defined above when $\fX = \Spf(\mathring\cA)$.
In particular, $\bar\fX$ is connected if $\fX^\an$ is connected.  The
inverse image of a closed point of $\bar\fX$ under the reduction map is
called a \emph{formal fiber}.

\paragraph[Shilov points]
\label{par:Shilovboundary}
The \emph{Shilov boundary} of a $K$-affinoid space $\sX = \sM(\cA)$
is defined to be the smallest closed subset $\Gamma(\sX) \subset\sX$ such that
every function $|f|$ for $f \in \cA$ attains its maximum at a point of
$\Gamma(\sX)$.  Let $\fX$ be a formal analytic variety over $\Spf(R)$
with Raynaud generic fiber $\sX$.  If $\eta$ is a generic point of
$\bar\fX$, then there is a unique preimage of $\eta$ under $\red$ which we
call the \emph{Shilov point} $x_{\eta}$ associated to $\eta$.
The residue field of $\eta$ is isomorphic to
$\td\sH(x_{\eta})$.  If $\sX = \sM(\cA)$ is affinoid then $\Gamma(\sX)$ is the
set of Shilov points of $\fX = \Spf(\mathring\cA)$.

\paragraph[Analytic curves] \label{par:analytic.curves}
Following \cite[\S{2.1.3}]{thuillier:thesis}, we define a (strictly)
\emph{analytic curve} over $K$ to be a (good) $K$-analytic space 
which is paracompact, of pure dimension
$1$, and without boundary.  The analytification of an algebraic curve over
$K$ (by which we mean a one-dimensional separated integral scheme of finite type over $K$)
is always an analytic curve in this sense.

If $\sX$ is an analytic curve and $\sV \subset\sX$ is an affinoid domain,
then by \cite[Proposition 3.1.3]{berkovich:analytic_geometry} and
\cite[Proposition 2.1.12]{thuillier:thesis} the following three (finite)
sets coincide: (i) the topological boundary $\del_\top\sV$ of $\sV$ in $\sX$; (ii) the
boundary $\partial\sV$ of $\sV$ in the sense of \cite[\S{2.5.7}]{berkovich:analytic_geometry};
and (iii) the Shilov boundary $\Gamma(\sV)$ of $\sV$.

\paragraph[Types of points in an analytic curve]
\label{par:types}
Let $x$ be a point in a $K$-analytic curve $\sX$ and let $\sH(x)$ be its
completed residue field.
The extension
$\td\sH(x) / k$ has transcendence degree $s(x) \leq 1$ and the
abelian group $|\sH(x)^\times|/|K^\times|$ has rank $t(x) \leq 1$.
Moreover, the integers
$s(x)$ and $t(x)$ must satisfy the \emph{Abhyankar inequality} 
\cite[Theorem~9.2]{vaquie:valuations}
\[
s(x) + t(x) \leq 1.
\]
Using the terminology from \cite{berkovich:analytic_geometry} and
\cite{berkovich:etalecohomology} (see also
\cite[\S{2.1}]{thuillier:thesis}), we say that $x$ is \emph{type-2} if
$s(x)=1$ and \emph{type-3} if $t(x)=1$.  If $s(x)=t(x)=0$, then $x$ is
called \emph{type-1} if $\sH(x) = K$
and \emph{type-4} otherwise.  Points of type 4 will not play any
significant role in this paper.  We define
\[\begin{split}
  \HH_\circ(\sX) &= \{\text{all points of } \sX \text{ of types } 2 
  \text{ and } 3 \} \\
  \HH(\sX) &= 
  \{\text{all points of } \sX \text{ of types } 2, 3, \text{ and } 4 \}.
\end{split}\]
We call $\HH_\circ(\sX)$ the set of \emph{skeletal points}, because it is the union of all skeleta of admissible formal models of $\fX$ 
(see \cite[Corollary 5.1]{bpr:berk_curves}),
$\HH(\sX)$ the set of \emph{norm points} of $\sX$, because it is the set of all points corresponding to norms on the function field $K(X)$ that extend the given norm on $K$.
If $\sX = X^{\an}$ is the analytification of an algebraic curve $X$ over $K$,
then $X(K) \subset X^{\an}$ is naturally identified with the set of type-1
points of $X^{\an}$, so $\HH(X^\an) = X^\an\setminus X(K)$. (Recall that
we are assuming throughout this discussion that $K$ is algebraically closed.)

\paragraph[Some facts about admissible $R$-algebras]
The following fact is standard and is easily proved using the results
of~\cite[\S1]{bl:fmI}. 

\begin{prop} \label{prop:compl.fin.pres}
  \begin{enumerate}
  \item If $A$ is a finitely presented and flat $R$-algebra then its
    $\varpi$-adic completion $\hat A$ is an admissible $R$-algebra.
  \item If $f:A\surject B$ is a surjective homomorphism of finitely
    presented and flat $R$-algebras with kernel $\fa$ then 
    $\hat f:\hat A\to\hat B$ is a surjection of admissible $R$-algebras
    with kernel $\fa\hat A$.
  \end{enumerate}
\end{prop}

\medskip

We set the following notation, which we will use
until~\parref{par:pure.degree}: $A$ and $B$ will denote admissible 
$R$-algebras, $\bar A = A\tensor_R k$ and $\bar B = B\tensor_R k$ their
reductions, and $A_K = A\tensor_R K$ and $B_K = B\tensor_R K$ the associated
$K$-affinoid algebras.  We let 
$\fX = \Spf(A)$, $\fY = \Spf(B)$, 
$\bar\fX=\Spec(\bar A)$, and $\bar\fY = \Spec(\bar B)$.  Let
$f: A\to B$ be a homomorphism, let $\bar f:\bar A\to\bar B$ and
$f_K: A_K\to B_K$ be the induced homomorphisms, and let
$\phi: \fY\to\fX$ and $\bar\phi:\bar\fY\to\bar\fX$ be the induced
morphisms. 

\begin{prop} \label{prop:finiteness.flatness}
  \begin{enumerate}
  \item $f$ is flat if and only if $\bar f: \bar A\to\bar B$ is flat.
  \item $f$ is finite if and only if $f_K: A_K\to B_K$ is finite.
  \end{enumerate}
\end{prop}

\pf 
The `only if' directions are clear.  Suppose that $\bar f$ is flat.
By~\cite[Lemma~1.6]{bl:fmI}, it suffices to show that
$f_n: A_n\to B_n$ is flat for all $n\geq 0$, where $A_n = A/\varpi^{n+1}A$
and $B_n = B/\varpi^{n+1}B$.  But $A_n$ and $B_n$ are of finite
presentation and flat over $R_n = R/\varpi^{n+1}R$, so $f_n$ is flat by
the fibral flatness criterion~\cite[Corollaire~11.3.1]{egaIV_3}.

Now suppose that $f_K$ is finite.  Choose a surjection 
$\mathring T_n\surject A$.  The induced homomorphism 
$T_n\surject A\tensor_R K\to B\tensor_R K$ is finite, so
by~\cite[Theorem~6.3.5/1]{bgr:nonarch} the composition
$\mathring T_n\surject A\to(B\tensor_R K)^\circ$ is integral.  Hence
$A\to(B\tensor_R K)^\circ$ is integral, so $A\to B$ is integral since 
$B\subset(B\tensor_R K)^\circ$.  Then
$f_n: A/\varpi^{n+1}A\to B/\varpi^{n+1} B$ is of finite type and integral for
all $n\geq 0$, so $f_n$ is finite, so $f$ is finite
by~\cite[Lemma~1.5]{bl:fmI}.\qed 

\begin{cor} \label{cor:finite.dominant}
  Suppose that $f_K:A_K\to B_K$ is finite and dominant, i.e.,\ that
  $\ker(f_K)$ is nilpotent.  Then $\bar f: \bar A\to\bar B$
  is finite and $\bar\phi:\bar\fY\to\bar\fX$ is surjective.
\end{cor}

\pf Since $A\subset A_K$ and $B\subset B_K$ we have that $\ker(f)$ is
nilpotent, and $f$ is finite by
Proposition~\ref{prop:finiteness.flatness}.  Hence $\Spec(B)\to\Spec(A)$
is surjective, so $\bar\fY\to\bar\fX$ is surjective.  Finiteness of $f$
implies finiteness of $\bar f$.\qed 

We say that a ring is \emph{equidimensional of dimension $d$}
provided that every maximal ideal  has height $d$.  Let
$\cA$ be a $K$-affinoid algebra, and let $\sX = \sM(\cA)$.  Then $\cA$ is
equidimensional of dimension $d$ if and only if $\dim(\sO_{\sX,x}) = d$ for
every $x\in\MaxSpec(\cA)$ by~\cite[Proposition~7.3.2/8]{bgr:nonarch}.
In particular, if $\sM(\cB)$ is an affinoid domain in $\sM(\cA)$ and $\cA$
is equidimensional of dimension $d$ then so is $\cB$.

\begin{prop} \label{prop:equidim.reduction}
  If $A_K$ is equidimensional of dimension $d$
  then $\bar A$ is equidimensional of dimension $d$.
\end{prop}

\pf Replacing $\fX$ with an irreducible formal affine open subset,
we may assume that $\fX$ is irreducible.  Let
$R\angles{x_1,\ldots,x_n}\surject A$ be a presentation of $A$.  By Noether
normalization~\cite[Theorem~6.1.2/1]{bgr:nonarch} we can choose the $x_i$
such that $K\angles{x_1,\ldots,x_d}\to A_K$ is finite and injective, where
$d = \dim(A_K)$.
Then $\bar\fX\to\A^d_k$ is finite and surjective
by Corollary~\ref{cor:finite.dominant}.\qed

\begin{cor} \label{cor:finite.surjective}
  Suppose that $f_K: A_K\to B_K$ is finite and dominant, and that $A_K$
  and $B_K$ are equidimensional (necessarily of the same dimension).  Then
  $\bar\phi:\bar\fY\to\bar\fX$ is finite and surjective, and
  the image of an irreducible component of $\bar\fY$ is an irreducible
  component of $\bar\fX$.
\end{cor}

\pf This follows immediately from Proposition~\ref{prop:equidim.reduction}
and Corollary~\ref{cor:finite.dominant}.\qed

The following theorem uses the fact that $K$ is algebraically closed in an
essential way.  It can be found 
in~\cite[Proposition~1.1]{bosch_lutk:uniformization}.

\begin{thm} \label{thm:admissible.reduced}
  Let $\cA$ be a $K$-affinoid algebra.  Then $\mathring\cA$ is admissible
  if and only if $\cA$ is reduced.
\end{thm}

\pf Since $\mathring\cA$ is always $R$-flat,
by~\cite[Proposition~1.1(c)]{bl:fmI} the issue is whether
$\mathring\cA$ is topologically finitely generated.  Suppose that
$\cA$ is reduced.  By~\cite[Theorem~6.4.3/1]{bgr:nonarch} there is a
surjection $T_n\surject\cA$ such that the residue norm on $\cA$ agrees
with the supremum norm; then by Proposition~6.4.3/3(i) of loc.\ cit.\ the
induced homomorphism $\mathring T_n\to\mathring\cA$ is surjective.
The converse follows in a similar way from Theorem 6.4.3/1 and Corollary~6.4.3/6 of loc.\
cit. \qed 

\begin{prop} \label{prop:canonical.reduction}
  The ring $\bar A$ is reduced if and only if $A = \mathring A_K$, in
  which case $A$ is reduced.
\end{prop}  

\pf 
If $A = \mathring A_K$ then $A$ is reduced by
Theorem~\ref{thm:admissible.reduced}, so
$\bar A = \td A_K$ is reduced.  Conversely,
suppose that $\bar A$ is reduced.
Let $\alpha: \mathring T_n\surject A$ be a surjection.
Since the $\mathring T_n$-ideal 
$\check T_n + \ker(\alpha) = \fm \mathring T_n + \ker(\alpha)$ is
the kernel of the composite homomorphism
$\mathring T_n \to A \to\bar A$, it is a reduced ideal; hence
by~\cite[Propositions~6.4.3/4,\,6.4.3/3(i)]{bgr:nonarch} we have
$A = \alpha(\mathring T_n) = \mathring A_K$.\qed

\begin{cor} If \label{cor:integral.reduction}
  $\bar A$ is an integral domain then $A_K$ is an integral domain and
  $|\cdot|_{\sup}$ is multiplicative.
\end{cor}

\pf By Proposition~\ref{prop:canonical.reduction} we have
$\td A_K = \bar A$, so the result follows
from \cite[Proposition 6.2.3/5]{bgr:nonarch}.\qed

\paragraphtoc[Finite morphisms of pure degree] \label{par:pure.degree}
In general there is not a good notion of the `degree' of a finite morphism
$Y\to X$ between noetherian schemes when $X$ is not irreducible, since the
degree of the induced map on an irreducible component of $X$ can vary from
component to component.   The notion of a
morphism having `pure degree' essentially means that the degree is the
same on every irreducible component of $X$.  This notion is quite well
behaved in that it respects analytification of algebraic varieties and of
admissible formal schemes.  The definition of a morphism of pure degree is
best formulated in the language of fundamental cycles.  We refer
to~\cite{thorup:rational_equivalence} for a review of the theory of cycles on a
noetherian scheme which is not necessarily of finite type over a field.

\paragraph
Let $X$ be a noetherian scheme.  A \emph{cycle} on $X$ is a finite
formal sum $\sum_W n_W\cdot W$, where $n_W\in\Z$ and $W$ ranges over
the irreducible closed subsets of $X$.  The group of cycles on $X$ is
denoted $C(X)$.  The \emph{fundamental cycle} of $X$ is
the cycle
\[ [X] = \sum_\zeta \length_{\sO_{X,\zeta}}(\sO_{X,\zeta}) \cdot \bar{\{\zeta\}}, \]
where the sum is taken over all generic points of $X$.  We define pushforwards
and pullbacks as in intersection theory.  These satisfy the usual properties: see~\cite[Lemmas~2.4, 2.5, and~4.8]{thorup:rational_equivalence}.

\begin{defn}
  Let $f: Y\to X$ be a finite morphism of noetherian schemes.  We say that
  $f$ has \emph{pure degree $\delta$} and we write
  $[Y:X] = \delta$ provided that $f_*[Y] = \delta\,[X]$; here
  $\delta\in\Q$ need not be an integer.  
\end{defn}

\begin{rem} \label{rem:pure.degree}
  Let $f: Y\to X$ be a finite morphism of noetherian schemes.
  \begin{enumerate}
  \item If $X$ is irreducible and every generic point of $Y$ maps to the
    generic point of $X$, then $f$ automatically has a pure degree, which
    we simply call the \emph{degree} of $f$.  Moreover, if $X$ is integral
    with generic point $\zeta$ then the degree of $f$ is the
    dimension of $\Gamma(f\inv(\zeta),\sO_{f\inv(\zeta)})$ as a vector
    space over the function field $\sO_{X,\zeta}$.  In particular, if $f$
    is a finite and dominant morphism of integral schemes, then the (pure)
    degree of $f$ is the degree of the extension of function fields.

  \item Let $\zeta$ be a generic point of $X$ and let
    $C=\bar{\{\zeta\}}$ be the corresponding irreducible component.
    Define the \emph{multiplicity of $C$ in $X$} to be the quantity
    \[ \mult_X(C) = \length_{\sO_{X,\zeta}}(\sO_{X,\zeta}), \]
    so $[X] = \sum_C \mult_X(C)\cdot C$.
    It follows that $f$ has pure degree $\delta$ if
    and only if (1) every irreducible component $D$ of $Y$ maps to an
    irreducible component of $X$, and (2) for every irreducible component $C$
    of $X$ we have
    \begin{equation} \label{eq:degree.alternate}      
      \delta\,\mult_X(C) = \sum_{D\surject C} \mult_Y(D)\,[D:C], 
    \end{equation}
    where $[D:C]$ is the usual degree of a finite morphism of integral
    schemes. 

  \item Let $g:X\to Z$ be another finite morphisms of noetherian
    schemes.  Suppose that $f$ has pure degree $\delta$ and $g$ has pure
    degree $\epsilon$.  Then $g\circ f$ has pure degree $\delta\epsilon$.

  \end{enumerate}
\end{rem}

\begin{prop} \label{prop:degree.basechange}
  Let $X,Y,X'$ be noetherian schemes,
  let $f:Y\to X$ be a finite morphism,
  let $g: X' \to X$ be a flat morphism, let 
  $Y' = Y\times_X X'$, and let $f':Y'\to X'$ be the projection.
  \begin{enumerate}
  \item If $f$ has pure degree $\delta$ 
    then $f'$ has pure degree $\delta$.
  \item If $g$ is surjective then $f$ has pure degree $\delta$ if and
    only if $f'$ has pure degree $\delta$.
  \end{enumerate}
\end{prop}

\pf Let $h:Y'\to Y$ be the other projection, so $h$ is flat.  
We have
\[ f'_*[Y'] = f'_*h^*[Y] = g^*f_*[Y] = \delta\, g^*[X] = \delta\,[X'], \]
which proves~(1).  Conversely, suppose that $g$ is surjective (and flat) and
that $f'_*[Y'] = \delta\,[X']$.  Then
\[ g^*f_*[Y] = f'_*h^*[Y] = f'_*[Y'] = \delta\,[X'] = g^*(\delta\,[X]), \]
so we are done because $g^*$ is visibly injective in this situation.\qed

\paragraph 
Next we will define pure-degree morphisms of analytic spaces.  As above,
we must first review the notion of the fundamental cycle of an analytic
space, as defined by Gubler~\cite[\S2]{gubler:local_heights}.

Let $\sX$ be a $K$-analytic space (assumed from now on to be Hausdorff and paracompact).
A \emph{Zariski-closed subspace} of
$\sX$ is by definition an isomorphism class of closed immersions
$\sV\inject\sX$.  A Zariski-closed subspace of $\sX$
is \emph{irreducible} if it cannot be expressed as a union of two proper
Zariski-closed subspaces.  Gubler~\cite[\S2]{gubler:local_heights} defines a
\emph{cycle} on $\sX$ to be a locally finite formal sum
$\sum_{\sV} n_\sV\,\sV$, where $n_\sV\in\Z$ and $\sV$ ranges over the
irreducible Zariski-closed subspaces of $\sX$; `locally finite' means that
there exists an admissible covering of $\sX$ by affinoid domains
intersecting only finitely many $\sV$ with $n_\sV\neq 0$.
Let $C(\sX)$ denote the group of cycles on $\sX$.  

\subparagraph
If $\sX=\sM(\cA)$ is affinoid then the Zariski-closed
subspaces of $\sX$ are in natural inclusion-reversing bijection with the
ideals of $\cA$; therefore we have an identification
$C(\sX) = C(\Spec(\cA))$, which we will make implicitly from now on.

\subparagraph
There are natural proper pushforward and flat pullback homomorphisms for cycles
on analytic spaces, which satisfy the expected properties.
There is a canonical
\emph{fundamental cycle} $[\sX]\in C(\sX)$ which is uniquely determined by
the property that for every affinoid domain 
$\iota: \sM(\cA)\inject\sX$, we have
$\iota^*[\sX]=[\sM(\cA)] = [\Spec(\cA)]$.  
See~\cite[2.6, 2.7, 2.8, and Proposition~2.12]{gubler:local_heights}.

\begin{defn}
  Let $f:\sY\to\sX$ be a finite morphism of $K$-analytic spaces.  We say
  that $f$ has \emph{pure degree $\delta$} and we write $[\sY:\sX]=\delta$
  provided that $f_*[\sY]=\delta[\sX]$.  Again $\delta\in\Q$ need not be
  an integer.
\end{defn}

\begin{rem} \label{rem:analytic.puredegree}
  Let $f:\sY\to\sX$ be a finite morphism of $K$-analytic spaces.
  \begin{enumerate}
  \item If $\sX=\sM(\cA)$ and $\sY=\sM(\cB)$ are affinoid then 
    $f:\sM(\cB)\to\sM(\cA)$ has
    pure degree $\delta$ if and only if the map of affine schemes
    $\Spec(\cB)\to\Spec(\cA)$ has pure degree $\delta$.

  \item If $f$ has pure degree $\delta$ and 
    $g:\sX\to\sZ$ is a finite morphism of analytic spaces of pure degree
    $\epsilon$ then $g\circ f$ has pure degree $\delta\epsilon$.
  \end{enumerate}
\end{rem}

\begin{prop} \label{prop:degree.cover}
  Let $f:\sY\to\sX$ be a finite morphism of $K$-analytic spaces.  
  \begin{enumerate}
  \item If $f$ has pure degree $\delta$, $\sM(\cA)\subset\sX$ is an
    affinoid domain, and $\sM(\cB) = f\inv(\sM(\cA))$, then 
    $\sM(\cB)\to\sM(\cA)$ has pure degree $\delta$.

  \item If there exists an admissible cover $\sX = \bigcup_i\sM(\cA_i)$
    of $\sX$ by affinoid domains such that 
    $\sM(\cB_i) = f\inv(\sM(\cA_i))\to\sM(\cA_i)$ has pure degree $\delta$
    for each $i$, then $f$ has pure degree $\delta$.
  \end{enumerate}
\end{prop}

\pf Since the inclusion $\sM(\cA)\inject\sX$ is flat, the first part
follows as in the proof of Proposition~\ref{prop:degree.basechange}(1).
In the situation of~(2), let
$f_i = f|_{\sM(\cB_i)}:\sM(\cB_i)\to\sM(\cA_i)$, and assume that 
$(f_i)_*[\sM(\cB_i)]=\delta[\sM(\cA_i)]$ for all $i$.  Arguing as in 
the proof of Proposition~\ref{prop:degree.basechange}(1), we see that
the pullback of $f_*[\sY]$ to $\sM(\cA_i)$ is equal to 
$\delta\,[\sM(\cA_i)]$ for all $i$; since $[\sX]$ is the unique cycle which
pulls back to $[\sM(\cA_i)]$ for all $i$, this shows that 
$f_*[\sY]=\delta\,[\sX]$.\qed

\medskip

The property of being a finite morphism of pure degree is compatible with
analytification: 

\begin{prop} \label{prop:degree.analytification}
  Let $f:Y\to X$ be a morphism of finite-type $K$-schemes.  
  Then $f$ is finite of pure degree $\delta$ if and only if
  $f^\an: Y^\an\to X^\an$ is finite of pure degree $\delta$.
\end{prop}

\pf By~\cite[Theorem~A.2.1]{conrad:irredcomps}, $f$ is finite if and only
if $f^\an$ is finite.  Hence we may assume that $X = \Spec(A)$ and 
$Y = \Spec(B)$ are affine.   If $\sM(\cA)\subset X^\an$ is an affinoid
domain then $\Spec(\cA)\to\Spec(A)$ is flat by Lemma~A.1.2 of loc.\ cit.,
and if $\{\sM(\cA_i)\}_{i\in I}$ is an admissible covering of $X^\an$ then 
$\Djunion_{i\in I}\Spec(\cA_i)\to\Spec(A)$ is flat and surjective.
Let $\sM(\cB_i)=f\inv(\sM(\cA_i))$.  We claim that
$\cB_i = B\tensor_A\cA_i$.  Since $B\tensor_A\cA_i$ is finite over $\cA_i$ it is
affinoid by~\cite[Proposition~6.1.1/6]{bgr:nonarch}, so the claim follows
easily from the universal property of the analytification (see
also~\cite[{\S}A.2]{conrad:irredcomps}).
Hence by Proposition~\ref{prop:degree.basechange}(2), 
$f$ has pure degree $\delta$ if and only if
$\Spec(\cB_i)\to\Spec(\cA_i)$ has pure degree $\delta$ for each $i$; 
by Remark~\ref{rem:analytic.puredegree}(1), this 
is the case if and only if 
$\sM(\cB_i)\to\sM(\cA_i)$ has pure degree $\delta$ for 
each $i$, which is equivalent to $f^\an$ having pure degree $\delta$ 
by Proposition~\ref{prop:degree.cover}(2).\qed

\medskip
The following counterpart to Proposition~\ref{prop:degree.analytification}
allows us to compare the degrees of the generic and special fibers of a
finite morphism of admissible formal schemes.  It will play a key role
throughout this paper.

\begin{prop}[Projection formula]
  \label{prop:projection.formula}
  Let $f:\fY\to\fX$ be a finite morphism of 
  admissible formal schemes, and 
  let $f^\an:\fY^\an\to\fX^\an$ and $\bar f:\bar\fY\to\bar\fX$ be the
  induced morphisms on the generic and special fibers, respectively.  If
  $f^\an$ has pure degree $\delta$ then $\bar f$ has pure degree
  $\delta$.
\end{prop}

\pf The theory of cycles on analytic spaces discussed above is
part of Gubler's more general intersection theory on admissible formal
schemes, and our `projection formula' is in fact a special case of Gubler's
projection formula~\cite[Proposition~4.5]{gubler:local_heights}; this can
be seen as follows.
Choose any $\varpi\in K^\times$ with $\val(\varpi)\in(0,\infty)$, and let $D$
be the Cartier divisor on $\fX$ defined by $\varpi$.  Essentially by
definition (cf.~(3.8) and~(3.10) of loc.\ cit.) the intersection
product $D.[\fX^\an]$ is equal to $\val(\varpi)\,[\bar\fX]$, and
likewise $(f^*D).[\fY^\an] = \val(\varpi)\,[\bar\fY]$.  Hence if 
$f^\an_*[\fY^\an] = \delta\,[\fX^\an]$ then
\[ \val(\varpi) \bar f_*[\bar\fY] = f_*((f^*D).[\fY^\an])
= D.f^\an_*[\fY^\an] = D.(\delta\,[\fX^\an]) 
= \val(\varpi)\,\delta\,[\bar\fX], \]
where the second equality is by Gubler's projection formula.
Canceling the factors of $\val(\varpi)$ yields
Proposition~\ref{prop:projection.formula}.\qed

\begin{rem}
  The converse to Proposition~\ref{prop:projection.formula} does not hold
  in general.  The following example is due to Gubler: let 
  $\fX = \Spf(R[x]/(x(x-\varpi)))$ and $\fY = \fX\djunion\Spf(R)$, and let 
  $f:\fY\to\fX$ be the map which is the identity on $\fX$ and which
  maps $\Spf(R)$ to $\fX$ via $x\mapsto 0$.  Then $f^\an$ does not have a
  pure degree, but $\bar f$ does since $\bar\fX$ is a point.
\end{rem}

\paragraph
Here we note some special cases of the projection formula:
\begin{enumerate}
\item Suppose that $\fX = \Spf(A)$ and $\fY = \Spf(B)$, and that $A$ is an
  integral domain with fraction field $Q$.
  If all generic points of $\Spec(B\tensor_R K)$ map to the generic point
  of $\Spec(A\tensor_R K)$ then
  $\sM(B\tensor_R K) \to \sM(A\tensor_R K)$ is finite with pure
  degree equal to $\dim_Q(B\tensor_A Q)$.
  By~\eqref{eq:degree.alternate}, for every irreducible
  component $\bar\fC$ of $\bar\fX$ we have
  \begin{equation} \label{eq:degree.special1}
    \dim_Q(B\tensor_A Q)\cdot \mult_{\bar\fX}(\bar\fC)
    = \sum_{\bar\fD\surject\bar\fC} \mult_{\bar\fY}(\bar\fD)\cdot[\bar\fD:\bar\fC],     
  \end{equation}
  where the sum is taken over all irreducible components $\bar\fD$ of $\bar\fY$
  that surject onto $\bar\fC$.
  
\item Suppose that $f^\an:\fY^\an\to\fX^\an$ is an isomorphism.  Then for every
  irreducible component $\fC$ of $\bar\fX$ we have
  \begin{equation} \label{eq:degree.special2}
    \mult_{\bar\fX}(\bar\fC) 
    = \sum_{\bar\fD\surject\bar\fC} \mult_{\bar\fY}(\bar\fD)\cdot[\bar\fD:\bar\fC], 
  \end{equation}
  where the sum is taken over all irreducible components $\bar\fD$ of $\bar\fY$
  that surject onto $\bar\fC$, because an isomorphism has pure degree 1.
\end{enumerate}

\section{Tropical integral models}
\label{section:tropicalintegralmodels}

We continue to assume that $K$ is an algebraically closed field which is
complete with respect to a nontrivial nonarchimedean valuation.

\begin{notn}
  Let $M\cong\Z^n$ be a lattice, with dual lattice $N = \Hom (M, \Z)$.  If
  $H$ is an additive subgroup of $\R$, we write $M_H$ for $M \tensor_\Z
  H$, so $N_H$ is naturally identified with $\Hom(M,H)$.  We write
  $\angles{\cdot,\cdot}$ to denote the canonical pairings $M \times N
  \rightarrow \Z$ and $M_\R \times N_\R \rightarrow \R$.
  
  Let $\T = \Spec K[M]$ be the torus over $K$ with character lattice $M$.
  For $u$ in $M$, we write $x^u$ for the corresponding character,
  considered as a function in $K[M]$.
\end{notn}

\paragraph[Extended tropicalization]
\label{par:extendedtrop}
A point $\| \cdot \|$ in $\T^\an$ naturally determines a real valued
linear function on the character lattice $M$, taking $u$ to 
$-\log\|x^u\|$.
The induced tropicalization map $\trop: \T^\an \rightarrow N_\R$ is
continuous, proper, and surjective \cite{payne:analytification}.  The
image of $\T(K)$ is exactly $N_{G}$, which is dense in $N_\R$ because $G$
is nontrivial and divisible.

More generally, if $\sigma$ is a pointed rational polyhedral cone in $N_\R$ and
$Y_\sigma = \Spec K[\sigma^\vee \cap M]$ is the associated affine toric
variety with dense torus $\T$, 
then there is a natural tropicalization map from $Y_\sigma$ to the space
of additive semigroup homomorphisms $\Hom(\sigma^\vee \cap M, \R \cup
\{\infty\})$ taking a point $\| \cdot \|$ to the semigroup map
$u\mapsto-\log \| x^u \|$, where $-\log(0)$ is defined to be $\infty$.
See \cite{payne:analytification, jdr:trop_ps} for further details.  We
write $N_\RR(\sigma)$ for the image of $Y_\sigma^\an$ under this
extended tropicalization map.

\begin{defn}
We say that a point in $N_\RR(\sigma)$ is \emph{$G$-rational} if it is in the
subspace $\Hom(\sigma^\vee \cap M, G \cup \{\infty\})$. 
\end{defn}
Note that the image of any $K$-rational point of $Y_\sigma$ is $G$-rational.

For any toric variety $Y_\Delta$, the tropicalization $N_\RR(\Delta)$ is
the union of the spaces $N_\RR(\sigma)$ for $\sigma$ in $\Delta$, glued along
the open inclusions $N_\RR(\tau) \subset N_\RR(\sigma)$ for $\tau
\preceq \sigma$.  The tropicalization maps on torus invariant affine opens
are compatible with this gluing, and together give a natural continuous,
proper, and surjective map of topological spaces $\trop:Y_\Delta^\an
\rightarrow N_\RR(\Delta)$.  Note that the vector space $N_\R$, which is
the tropicalization of the dense torus $\T \subset Y_\Delta$, is open and
dense in $N_\RR(\Delta)$.  For the purpose of constructing tropical
integral models of toric varieties and their subvarieties, it will
generally suffice to study polyhedral complexes in $N_\R$.

Let $X$ be a closed subscheme of $Y_\Delta$.  The tropicalization
$\Trop(X)$ is the image of $X^\an$ under $\trop$.  Since $X(K)$ is dense
in $X^\an$, its image is dense in $\Trop(X)$.  Furthermore, every
$G$-rational point of $\Trop(X)$ is the image of a point of $X(K)$, and if
$X$ is irreducible then 
the preimage of any point in $\Trop(X) \cap N_{G}$ is Zariski dense in
$X$.  See~\cite[Corollary~4.2]{payne:fibers} and
\cite[Remark~2]{payne:fibers_correction},
\cite[Proposition~4.14]{gubler:tropical_guide},
or~\cite[Theorem~4.2.5]{osserman_payne:lifting}.

\paragraph[Polyhedral domains] \label{par:polytopal.domains} 
Recall that the \emph{recession cone} $\sigma_P$ of a nonempty polyhedron $P \subset N_\R$
is the set of those $v$ in $N_\R$ such that $w + v$ is in $P$ whenever $w$ is in $P$.
If $P$ is the intersection of the halfspaces $\angles{u_1,v} \geq a_1,
\ldots, \angles{u_r,v} \geq a_r$ then $\sigma_P$ is the dual of the cone
in $M_\R$ spanned by $u_1, \ldots, u_r$.  In particular, if $P$ is an
integral $G$-affine polyhedron, then these halfspaces can be chosen with
each $u_i$ in $M$, so the recession cone $\sigma_P$ is a rational polyhedral cone.
The recession cone can also be characterized as the intersection with $N_\R
\times \{0\}$ of the closure in $N_\R \times \R$ of the cone spanned by
$P \times \{1\}$. 

Let $P$ be an integral $G$-affine polyhedron in $N_\R$ that does not contain any positive dimensional affine linear subspace, 
so its recession cone $\sigma = \sigma_P$ is pointed.

\begin{defn}
  The \emph{polyhedral domain} associated to $P$ is the inverse image under 
  $\trop: Y_\sigma^\an\to N_\RR(\sigma)$ of the
  closure of $P$ in $N_\RR(\sigma)$ and is denoted $\sU^P$.
\end{defn}

These polyhedral domains, introduced in~\cite{jdr:trop_ps},
directly generalize the polytopal domains studied
by Gubler in \cite{gubler:tropical}.  Indeed, a polytopal domain is the
preimage in $\T^\an$ of an integral $G$-affine polytope in $N_\R$.  Since
the recession cone of a polytope in $N_\R$ is the zero cone, whose
associated toric variety is $\T$, Gubler's polytopal domains are exactly
the special case of these polyhedral domains where $P$ is bounded. 

By~\cite[\S6]{jdr:trop_ps} the polyhedral domain $\sU^P$ is an affinoid
domain in $Y_\sigma^\an$ with coordinate ring 
\[ 
K\angles{\sU^P} = \left\{\sum_{u\in \sigma^\vee \cap M} a_u x^u~:~ \lim(\val(a_u) + \angles{u,v}) = \infty \text{ for all } v\in P\right\}, 
\]
where the limit is taken over all complements of finite sets.  Its supremum norm is given by
\begin{equation} \label{eq:sup.norm}
\left|\sum a_u x^u\right|_{\sup} = 
\sup_{\substack{u \in\sigma^\vee\cap M\\v \in P}} |a_u|\exp(-\angles{u,v}).
\end{equation}
Since the recession cone $\sigma$ is pointed, the polyhedron $P$ contains no linear subspace and hence has vertices.  
The supremum above is always achieved at one of the vertices of $P$, so the ring of power-bounded regular functions on $\sU^P$ is
\begin{equation} \label{eq:power.bounded}
K\angles{\sU^P}^\circ = 
\left\{\sum_{u\in \sigma^\vee\cap M} a_u x^u\in K\angles{\sU^P} ~:~ \val(a_u) + \angles{u,v} \geq 0
\text{ for all } v\in\vertices(P)\right\}. 
\end{equation}
Since $K \angles{\sU^P}$ is reduced, Theorem~\ref{thm:admissible.reduced}
implies that $\fU^P = \Spf(K\angles{\sU^P}^\circ)$ is an admissible formal
scheme with analytic generic fiber $\sU^P$.  

\begin{rem} \label{rem:non.strict.polydomains}
  If $P$ is integral affine but not $G$-affine then the inverse image
  $\sU^P$ of the closure of $P$ under $\trop$ is a non-strict affinoid
  domain.  Indeed, if $K'$ is a complete valued field extension of $K$ whose value
  group $G'$ is large enough that $P$ is $G'$-affine then
  $\sU^P\hat\tensor_K K'$ is strictly $K'$-affinoid.
\end{rem}

\paragraph[Polyhedral integral models]
Let $P$ be an integral $G$-affine polyhedron in $N_\R$ whose recession
cone $\sigma = \sigma_P$ is pointed.
As usual, we let $Y_\sigma = \Spec K[\sigma^\vee \cap M]$ denote
the associated affine toric variety with dense torus $\T$.

\begin{defn}
  We define $R[Y^P] \subset K[\sigma^\vee \cap M]$ to be the subring
  consisting of those Laurent polynomials $\sum a_u x^u$ such that
  $\val(a_u) + \angles{u,v}\geq 0$ for all $v \in P$ and all $u$.
  The scheme $Y^P \coloneq \Spec(R[Y^P])$ is called a 
  \emph{polyhedral integral model} of $Y_\sigma$.
\end{defn}

In other words, $R[Y^P]$ is the intersection of $K\angles{\sU^P}^\circ$
with $K[M]$.  It is clear that $K\angles{\sU^P}^\circ$ is the
$\varpi$-adic completion of $R[Y^P]$.  Note that $R[Y^P]$ is torsion-free
and hence flat over $R$.

\begin{lem}
The tensor product $R[Y^P] \tensor_R K$ is equal to $K[Y_\sigma]$.
\end{lem}

\proof By definition we have $R[Y^P]\tensor_R K\subset K[Y_\sigma]$.  For
the other inclusion, note that if $g = \sum b_u x^u$ is in $K[Y_\sigma]$
then the minimum over $v$ in $P$ of $\val(b_u) + \angles{u,v}$ is achieved
at some vertex of $P$.  It follows that some sufficiently high power of
$\varpi$ times $g$ is in $R[Y^P]$, and hence $g$ is in $R[Y^P] \tensor_R
K$.  \qed 

\begin{rem}
One could equivalently define $R[Y^P]$ to be the subring of $K[M]$ satisfying the same inequalities.  
Since $P$ is closed under addition of points in $\sigma$, any Laurent polynomial satisfying these inequalities for all $v$ in $P$ 
must be supported in $\sigma^\vee$.
\end{rem}

\medskip

We will use the following notation in the proof of
Proposition~\ref{prop:fg} below.
For each face $F \leq P$, let
$\sigma(F)$ be the cone in $N_\R$ spanned by $P-v$ for any $v$ in the
relative interior of $F$.  In other words, $\sigma(F) = \Star_P(F)$.  We
fix a labeling $v_1, \ldots, v_r$ for the vertices of $P$, and write
$\sigma_i$ for $\sigma(v_i)$.  The dual cone $\sigma_i^\vee$ is
\[
\sigma_i^\vee = \{ u \in \sigma_P^\vee \ : \ \angles{u, v_i} \leq \angles {u, v_j} \text{ for all } j \}.
\]
The cones $\sigma_1^\vee, \ldots, \sigma_r^\vee$ are the maximal cones of the (possibly degenerate) inner normal fan of $P$, 
and their union is $\sigma_P^\vee$.

\begin{prop} \label{prop:fg}
Let $P$ be a $G$-rational polyhedron in $N_\R$. Then $R[Y^P]$ is finitely presented over $R$.
\end{prop}

\pf
By~\cite[Corollary~3.4.7]{raynaud_gruson:platitude}, any finitely generated and flat algebra over an integral domain 
is automatically of finite presentation, so it suffices to show that $R[Y^P]$ is finitely generated. 

The cones $\sigma_1^\vee, \ldots, \sigma_r^\vee$ cover $\sigma^\vee$, so $R[Y^P]$ is generated by the subrings
\[
A_j = R[Y^P] \cap K[\sigma_j^\vee \cap M]
\]
for $1 \leq j \leq r$.  Therefore, it will suffice to show that each $A_j$ is finitely generated over $R$.

The semigroup $\sigma_j^\vee \cap M$ is finitely generated by Gordan's Lemma  \cite[p.~12]{fulton:toric_varieties}.  
Let $u_1, \ldots, u_s$ be generators, and choose $a_1, \ldots, a_s$ in $K^\times$ such that $\val(a_i) + \angles{u_i, v_j} = 0$.  
Then each monomial in $A_j$ can be written as an element of $R$ times a monomial in the $a_i x^{u_i}$.  
It follows that $A_j$ is finitely generated over $R$, as required, with generating set $\{a_1 x^{u_1}, \ldots, a_s x^{u_s}\}$.
\qed

In particular, $Y^P$ is a flat and finitely presented $R$-model of
the affine toric variety $Y_\sigma$.

\begin{rem}
  As in Remark~\ref{rem:non.strict.polydomains}, one can construct
  an algebraic model $Y^P$ of $Y_\sigma$ associated to an integral affine
  but not $G$-affine polyhedron $P$; when $P$ is a point this is done
  in~\cite{osserman_payne:lifting}.  This model is not of finite type.
\end{rem}

\paragraph[Polyhedral integral and formal models of subschemes]
\label{par:integral.models}
Let $P$ be an integral $G$-affine polyhedron with pointed recession cone
$\sigma$.  Let $X$ be the closed subscheme of the affine toric variety
$Y_\sigma$ over $K$ defined by an ideal $\fa\subset K[Y_\sigma]$.

\begin{defn} \label{defn:polyhedral.models}
  \begin{enumerate}
  \item Let $\sX^P = X^\an\cap\sU^P$.  This is the Zariski-closed subspace of
    $\sU^P$ defined by $\fa K\angles{\sU^P}$.
  \item The \emph{polyhedral integral model} of $X$ is the
    scheme-theoretic closure $X^P$ of $X$ in $Y^P$.  
    It is defined by the ideal $\fa^P = \fa \cap R[Y^P]$. 
  \item The \emph{polyhedral formal model} of $\sX^P$ is the $\varpi$-adic
    completion $\fX^P$ of $X^P$.  We will show
    in Proposition~\ref{prop:fXbar.is.inPX} that $\fX^P$ is an admissible formal
    scheme with generic fiber $\sX^P$.
  \item The \emph{canonical model} of $\sX^P$ is 
    \[ \fX^P_\can = \Spf\big((K\angles{\sU^P}/\fa
    K\angles{\sU^P})^\circ\big). \]
    By Theorem~\ref{thm:admissible.reduced}, the canonical model is admissible if and only if $\sX^P$ is reduced.  
  \end{enumerate}
\end{defn}

\begin{notn}
  The \emph{$P$-initial degeneration} of $X$ is defined to be
  \[ \inn_P(X) = X^P\tensor_R k = \fX^P\tensor_R k. \]
\end{notn}

As usual we write $\bar\fX{}^P_\can = \fX^P_\can\tensor_R k$.  This coincides
with the canonical reduction of $\sX^P$ when $\sX^P$ is reduced.
In the case where $P$ is a single point $w\in N_G$ we write
$\sX^w$, $\fX^w$, $\inn_w(X)$, etc.  In this case, $\inn_w(X)$ is the
$w$-initial degeneration of $X$ in the sense generally used in the
literature (and in the introduction). 

\begin{lem}
  The ideal $\fa^P$ is finitely generated.
\end{lem}

\pf Since $X^P$ is the closure of its generic fiber, it is flat over
$\Spec R$, and its coordinate ring is a quotient of the finitely generated
$R$-algebra $R[M]$.  Since any finitely generated flat algebra over an
integral domain is finitely presented
\cite[Corollary~3.4.7]{raynaud_gruson:platitude}, it follows that $\fa^P$
is finitely generated.  \qed

\begin{prop} \label{prop:fXbar.is.inPX}
  The formal scheme $\fX^P$ is the formal closed subscheme of 
  $\fU^P$ defined by $\fa^P K\angles{\sU^P}^\circ$.  It is an admissible
  formal scheme with generic fiber 
  $\sX^P$ and special fiber $\inn_P(X)$.
\end{prop}

\pf The admissibility of $\fX^P$ is a consequence
of Proposition~\ref{prop:compl.fin.pres}(1). 
If $A = R[Y^P]/\fa$ then by definition $X^P = \Spec(A)$ and
$\fX^P = \Spf(\hat A)$, where $\hat A$ is the $\varpi$-adic completion of $A$.
By Proposition~\ref{prop:compl.fin.pres}(2) the sequence
\[ 0 \To \fa^P K\angles{\sU^P}^\circ \To K\angles{\sU^P}^\circ\To\hat A\To
0 \]
is exact; it follows that $\fX^P$ is the closed subscheme of $\fU^P$
defined by $\fa^P K\angles{\sU^P}^\circ$.  We have
\[ \big(K\angles{\sU^P}^\circ / \fa^P K\angles{\sU^P}^\circ\big)\tensor_R K
= K\angles{\sU^P}/\fa K\angles{\sU^P} \]
since $K\fa^P = \fa$, so $\fX^P\tensor_R K = \sX^P$.  The special fiber of
$\fX^P$ agrees with the special fiber of $X^P$ by construction.\qed

The canonical inclusion
\[ K\angles{\sU^P}^\circ/\fa^P K\angles{\sU^P}^\circ \inject
\big(K\angles{\sU^P}/\fa K\angles{\sU^P}\big)^\circ \]
induces a map of formal schemes
\[ \fX^P_\can \To \fX^P. \]
As the above morphism induces an isomorphism on analytic generic fibers,
it is \emph{finite} when $\sX^P$ is reduced by
Proposition~\ref{prop:finiteness.flatness}(2). 
The special fiber of the above morphism is a morphism
$\bar\fX{}^P_\can\to\inn_P(X)$.
Many of the results of this paper are proved by using this morphism
and the results of \S\ref{section:admissible.R.algebras} (in particular
the projection formula, Proposition~\ref{prop:projection.formula}) to compare these two models.

\paragraph[Compatibility with extension of the ground field]
\label{par:base.change}
We continue to use the notation of~\parref{par:integral.models}.
Let $K'$ be an algebraically closed complete valued field extension of
$K$, with valuation ring $R'$ and residue field $k'$.  Let
$P$ be an integral $G$-affine polyhedron in $N_\R$ with pointed recession
cone $\sigma$.  Let $Y_\sigma' = Y_\sigma\tensor_K K'$, so $Y_\sigma'$ is
the affine toric variety defined over $K'$ with dense torus
$\T'\coloneq\T\tensor_K K'$ associated to the cone $\sigma$.
The triangle
\[\xymatrix @=.2in{
  {(Y_\sigma')^\an} \ar[rr] \ar[dr]_\trop & & {Y_\sigma^\an} \ar[dl]^\trop \\
  & {N_\R(\sigma)}
}\]
commutes, so $\sU^P\hat\tensor_K K'$ is the polyhedral domain in 
$(Y_\sigma')^\an$ associated to $P$.  Likewise the
polyhedral integral model $(Y')^P$ of $Y_\sigma'$ associated to $P$ is
naturally identified with $Y^P\tensor_R R'$.  Indeed, as an $R$-module we
have
\[ R[Y^P] = \Dsum_{u\in\sigma^\vee\cap M} 
R_u\cdot x^u \subset K[\sigma^\vee\cap M] \sptxt{where}
R_u = \Big\{ a\in R~:~\val(a)\geq\max_{v\in\vertices(P)}-\angles{u,v} \Big\}. \] 
Since $\angles{u,v}\in G$ for all $u\in M$ and $v\in\vertices(P)$ each 
$R_u$ is a free $R$-module of rank $1$, so the image of $R_u\tensor_R R'$ in 
$K'$ is exactly $R'_u$.

Let $X\subset Y_\sigma$ be the closed subscheme defined by an ideal
$\fa\subset K[\sigma^\vee\cap M]$ and 
let $X' = X\tensor_K K'\subset Y_\sigma'$,
so $X'$ is defined by $\fa K'[\sigma^\vee\cap M]$.  Since the above
triangle is 
commutative, we have $\Trop(X) = \Trop(X')\subset N_\R(\sigma)$, and
$\trop:(X')^\an\to\Trop(X)$ factors through the natural map
$(X')^\an\to X^\an$.  Hence
\[ (\sX')^P = \trop\inv(P)\cap(X')^\an = \sX^P\hat\tensor_K K'. \]
Since schematic closure commutes with flat base change, the polyhedral
integral model $(X')^P$ of $X'$ coincides with $X^P\tensor_R R'$; hence if
$\fa^P = \fa\cap R[Y^P]$ is the ideal defining $X^P$ then $(X')^P$ is
defined by $\fa^P R'[(Y')^P]$.  It follows from this and
Proposition~\ref{prop:fXbar.is.inPX} that 
$(\fX')^P = \fX^P\hat\tensor_R R'$, and in particular that 
$\inn_P(X') = \inn_P(X)\tensor_k k'$.  As for the canonical models,
suppose that $X$ is reduced, so $X'$ is reduced as well.  Then
$(\fX')^P_\can = \fX^P_\can\hat\tensor_R R'$ because
$(\fX^P_\can\hat\tensor_R R')\tensor_{R'} k'
= (\fX^P_\can\tensor_R k)\tensor_k k'$
is reduced; cf.\ Proposition~\ref{prop:canonical.reduction}.

Below we will make various definitions by passing to a valued field
extension $K'$ of $K$.  In order for these definitions to be independent
of the choice of $K'$, we will need the following fact, proven 
in~\cite[0.3.2]{ducros:excellents} or~\cite[\S4]{conrad:aws}.

\begin{lem} \label{lem:dominating.extension}
  Let $K_1,K_2$ be complete valued field extensions of $K$.  Then there is
  a complete valued field extension $K'$ of $K$ admitting isometric embeddings 
  $K_1\inject K'$ and $K_2\inject K'$ over $K$.
\end{lem}

\paragraphtoc[Relative multiplicities and tropical multiplicities]
\label{par:relative.mult}
Recall~\parref{par:Shilovboundary} that if $\sX = \sM(\cA)$ is an affinoid
space then the reduction map induces a one-to-one correspondence between
the Shilov boundary points of $\sX^\an$ and the generic points of the
canonical reduction $\Spec(\td\cA)$.  This leads to the following
definition:

\begin{defn} \label{defn:m.rel.xi}
  Let $X\subset\T$ be a reduced and equidimensional closed subscheme, let 
  $x\in X^\an$, let $w=\trop(x)$, and suppose that $w\in N_G$.  Define
  the \emph{relative multiplicity $m_{\rel}(x)$ of $x$ in $\trop\inv(w)$}
  as follows.  If $x$ is not a Shilov boundary point of $\trop\inv(w)$
  then we define its multiplicity to be zero.  Otherwise $\red(x)$ is
  the generic point of an irreducible component $\bar\fC$ of 
  $\bar\fX{}^w_\can$; we define the multiplicity of $x$ to be
  $[\bar\fC:\im(\bar\fC)]$, where $\im(\bar\fC)$ 
  is the image of $\bar\fC$ in $\inn_w(X)$ (this is an irreducible component
  by Corollary~\ref{cor:finite.surjective}).

  Now suppose that $w\notin N_G$.  Let $K'$ be an algebraically closed
  complete valued field extension of $K$ such that $w\in N_{G'}$, where
  $G'$ is the value group of $K'$.  Let $X' = X\tensor_K K'$ 
  and let $\phi:(X')^\an\to X^\an$ be the natural morphism.  We define
  \[ m_\rel(x) = \sum_{x'\in\phi\inv(x)} m_\rel(x'). \]
\end{defn}

In order for the above definition to make sense, by
Lemma~\ref{lem:dominating.extension} we only have to show that if 
$K\subset K'\subset K''$ are algebraically closed complete valued field
extensions then we can calculate $m_\rel(x)$ with respect to either $K'$
or $K''$.  Replacing $K$ with $K'$, we are reduced to showing:

\begin{lem} \label{lem:multiplicities.extension.scalars}
  Let $K'$ be an algebraically closed complete
  valued field extension of $K$ and let $X' = X\tensor_K K'$.  Let $x\in
  X^\an$, let  $w = \trop(x)$, and suppose that $w\in N_G$.  Then the
  natural map $(\sX')^w\to\sX^w$ induces a bijection of Shilov boundary points
  which preserves relative multiplicities.
\end{lem}

\pf Let $k'$ be the residue field of $K'$.  As discussed
in~\parref{par:base.change} we have 
$\inn_w(X') = \inn_w(X)\tensor_k k'$ and
$(\bar\fX\p)^w_\can = \bar\fX\s{w}_\can\tensor_k k'$, so the first assertion
follows from the fact that 
$\bar\fX\s{w}_\can\tensor_k k'\to\bar\fX\s{w}_\can$ 
induces a bijection on irreducible components.
Let $\bar\fC$ be an irreducible component of $\bar\fX\s{w}_\can$ and let
$\bar\fD$ be its image in $\bar\fX\s{w}$.  Then 
$[\bar\fC:\bar\fD] = [\bar\fC\tensor_k k':\bar\fD\tensor_k k']$, so
relative multiplicities are preserved as well.\qed

Later we will relate $m_\rel(x)$ to other geometrically-defined notions
of multiplicity; see Proposition~\ref{prop:deg.from.mrel} and
Theorem~\ref{thm:compatible.multiplicities}.  For the moment we relate
relative multiplicities to tropical multiplicities, defined as follows:

\begin{defn}
  Let $X\subset\T$ be a closed subscheme and let $w\in\Trop(X)$.
  If $w\in N_G$ then the 
  \emph{tropical multiplicity of $X$ at $w$} is defined to be 
  \[ m_{\Trop}(w) = \sum_{\bar C\subset\inn_w(X)} \mult_{\inn_w(X)}(\bar C), \]
  where the sum is taken over all irreducible components $\bar C$ 
  of $\inn_w(X)$.  If $w\notin N_G$ then let $K'$ be an algebraically
  closed complete valued field extension of $K$ such that
  $w\in N_{G'}$, where $G'$ is the value group of $K'$.  Let 
  $X' = X\tensor_K K'$.  We define $m_{\Trop}(w)$ to be the tropical
  multiplicity of $w$ relative to 
  $\trop: (X')^\an\to\Trop(X)$.
\end{defn}

The fact that $m_{\Trop}(w)$ is independent of the choice of $K'$ is
proved in~\cite[Remark~A.5]{osserman_payne:lifting}.  It is also one
of the consequences of the following proposition.

\begin{prop}
  \label{prop:relativemults}
  Let $X\subset\T$ be a reduced and equidimensional closed subscheme and let 
  $w\in\Trop(X)$.  Then
  \[ m_{\Trop}(w) = \sum_{x\in\trop\inv(w)} m_{\rel}(x). \]
\end{prop}

\pf We immediately reduce to the
case where $w\in N_G$ by extending the ground field if necessary.   By
definition we have 
\[ \sum_{x\in\trop\inv(w)} m_{\rel}(x)
= \sum_{\bar\fC\subset\bar\fX{}^w_\can} [\bar\fC:\im(\bar\fC)] \]
where the sum is taken over all irreducible components $\bar\fC$ of
$\bar\fX{}^w_\can$; the image $\im(\bar\fC)$ of $\bar\fC$ in $\inn_w(X)$
is an irreducible component by Corollary~\ref{cor:finite.surjective}.
Also by definition,
\[ m_{\Trop}(w) = \sum_{\bar C\subset\inn_w(X)} \mult_{\inn_w(X)}(\bar C) \]
where the sum is taken over all irreducible components $\bar C$ of
$\inn_w(X)$.  By the projection
formula~\parref{prop:projection.formula}, for every 
irreducible component $\bar C$ of $\inn_w(X)$ we have
\[ \mult_{\inn_w(X)}(\bar C) = \sum_{\bar\fC\surject\bar C} [\bar\fC:\bar C] \]
where the sum is taken over all irreducible components $\bar\fC$ of 
$\bar\fX{}^w_\can$ mapping onto $\bar C$ (for any such $\bar\fC$ we have
$\mult_{\bar\fX{}^w_\can}(\bar\fC)=1$ since $\bar\fX{}^w_\can$ is
reduced).  Therefore 
\[ m_{\Trop}(w) = \sum_{\bar C\subset\inn_w(X)} \mult_{\inn_w(X)}(\bar C)
= \sum_{\bar C\subset\inn_w(X)}\sum_{\bar\fC\surject \bar C} [\bar\fC:\bar
C]
= \sum_{\bar\fC\subset\bar\fX{}^w_\can} [\bar\fC:\im(\bar\fC)]. \]
\qed

\begin{eg}
In Example \ref{eg:Example2b}, 
we have $m_{\rel}(\xi) = 1$ for all $\xi \in \Gamma$ and
$m_{\rel}(\xi) = 0$ for all $\xi \not\in \Gamma$.
This follows from Proposition~\ref{prop:relativemults} and the concrete
description of the tropicalization map in Example~\ref{eg:Example2b},
together with the observation that $m_{\Trop}(0,0)=m_{\Trop}(2,2) = 1$  
(since the initial degenerations 
$\inn_{(0,0)}(E) \cong \Spec k[x,y,x^{-1},y^{-1}]/(y^2 - x^3 - x^2)$ and $\inn_{(2,2)}(E) \cong \Spec k[x,y,x^{-1},y^{-1}]/(y^2 - x^2 - 1)$ 
are both integral schemes over $k$).  Note that $m_{\rel}(\xi)>0$ for all
$\xi\in\Gamma$ because $\xi$ is contained in the topological boundary of
$\trop\inv(\trop(\xi))$ in $E^\an$, hence in the Shilov boundary;
see~\parref{par:analytic.curves}. 

\end{eg}

\paragraph[Polyhedral structures on tropicalizations]
\label{par:polyhedral.structures}
Let $W$ be a $G$-rational affine space in $N_\R$ and let $W_0$ be the
linear space under $W$, so $W_0$ is spanned 
by $W_0\cap N$.  Set $N' = N/(W_0\cap N)$ and 
$M' = W_0^\perp\cap M\subset M$, and let $\T'$ be the torus 
$\Spec(K[M'])$.  We call $\T'$ the \emph{torus transverse to $W$}.
Let $w'\in N'_G$ be the image of any point of $W$.  Then
\[ R[(\T')^{w'}] = \bigg\{\sum_{u\in M'} a_u x^u\in K[M']~:~
\val(a_u) + \angles{u,w'}\geq 0\bigg\}, \]
so for all $w\in N_G\cap W$ we have $R[(\T')^{w'}]\subset R[\T^w]$.
Hence we have a natural morphism
$\pi_w:\bar\T{}^w\to(\bar\T{}')^{w'}$ for all $w\in N_G\cap W$.

\begin{rem} \label{rem:split.torus}
  Let $N'' = \ker(N\to N') = W_0\cap N$ and let 
  $M'' = \Hom_\Z(N'',\Z)$, so we have exact sequences
  \[ 0\To N''\To N\To N'\To 0 \sptxt{and}
  0\To M'\To M \To M''\To 0 \]
  inducing an exact sequence of tori
  \[ 0\To \T''\To\T\To\T'\To 0 \]
  where $\T'' = \Spec(K[M''])$.  We call $\T''$ the
  \emph{torus parallel to $W$}.
  Choosing a splitting of $N\surject N'$ splits all three exact sequences,
  and in particular furnishes an isomorphism $\T\cong\T'\times\T''$.
  Let $w\in W\cap N_G$ and let $w''$ be its image in $N''_\R$.
  Then we have an isomorphism
  $\bar\T{}^w\cong(\bar\T{}')^{w'}\times(\bar\T{}'')^{w''}$ under which
  $\pi_w$ corresponds to the projection onto the first
  factor.
\end{rem}

\begin{thm} \label{thm:polyhedral.decomp}
  Let $X\subset\T$ be an equidimensional subscheme of dimension $d$.  The
  set $\Trop(X)$ admits a polyhedral complex structure of pure dimension
  $d$ with the following properties:
  \begin{enumerate}
  \item The tropical multiplicities are constant along the relative
    interior of every maximal face.
  \item Let $w$ be contained in the relative interior of a maximal face
    $\tau$ of $\Trop(X)$, let $W = \spn(\tau)$, let $\T'$ be the torus
    transverse to $W$, and let
    $\pi_w:\bar\T{}^w\to(\bar\T{}')^{w'}$ be the natural map.  Then
    $\inn_w(X)\cong\pi_w\inv(\bar Y)$ for some dimension-zero
    subscheme $\bar Y$ of $(\bar\T{}')^{w'}$.
  \end{enumerate}
\end{thm}

\pf The first part is a basic result in tropical geometry; it is
proved in \cite[\S{3.3}]{maclagan_sturmfels:book}.%
\footnote{The proofs in \cite[\S{3.3}]{maclagan_sturmfels:book} assume
  that there is a section to the valuation map $\val : K^\times \to G$.
  Such a section always exists when $K$ is an algebraically closed
  nonarchimedean field; the following short proof was communicated to us
  by David Speyer.  If $G=\{ 0 \}$ then there is nothing to prove.
  Otherwise, consider the short exact sequence $0 \to U \to K^\times \to G
  \to 0$.  Since $K$ is algebraically closed, the group $U$ is
  divisible. Thus $U$ is injective as a $\ZZ$-module, so ${\rm Ext}^1(A,
  U)=0$ and the valuation map splits.}  
Let $\T''\subset\T$ be the torus parallel to $W$ and write
$\bar\T\pp = (\bar\T\pp)^0$, so 
$(\bar\T{}')^{w'}$ is the quotient of $\bar\T{}^w$ by $\bar\T\pp$.
By~\cite[Proposition~2.2.4]{speyer:thesis}, the initial degeneration
$\inn_w(X)$ is invariant under the action of $\bar\T\pp$.  Therefore
$\inn_w(X)$ is the inverse image of a closed subscheme $\bar Y$ of 
$(\bar\T\p)^{w'}$ (in fact $\bar Y$ is the quotient
$\inn_w(X)/\bar\T\pp\subset(\bar\T\p)^{w'}$); counting dimensions, we see
that $\dim(\bar Y) = 0$.  
See~\cite[Theorem~1.1.1 and Amplification~1.1.3]{mumford:git} for basic
existence results about geometric quotients of affine schemes over a field
by a free action of a reductive group.\qed

In the situation of Theorem~\ref{thm:polyhedral.decomp}(2),
let $\T''$ be the torus parallel to $W$, and choose a splitting
$\T\cong\T'\times\T''$ as in Remark~\ref{rem:split.torus}.
Then $\inn_w(X)\cong \bar Y\times (\bar\T{}'')^{w''}$.

\paragraphtoc[The tropical projection formula]
\label{par:finiteflatfibration}
Let $X\subset\T$ be a reduced and equidimensional closed subscheme of
dimension $d$ and let $P$ be an integral $G$-affine polytope contained in
the relative interior of a maximal ($d$-dimensional) face $\tau$ of a
polyhedral complex decomposition of $\Trop(X)$ as in
Theorem~\ref{thm:polyhedral.decomp}.  Let $W$ be the affine span of
$\tau$, let $\T'$ be the torus transverse to
$W$~\parref{par:polyhedral.structures}, let $\T''$ be the torus
parallel to $W$ (Remark~\ref{rem:split.torus}), and choose a splitting
$\T\to\T''$.
Note that $\dim(\T'') = d$.  Let $P''$ be the image
of $P$ in $N''_\R$, so $\sU^{P''}$ is a polytopal domain in $(\T'')^\an$.
The map $\sU^P\to\sU^{P''}$ induces a morphism $\psi^P:
\sX^P\to\sU^{P''}$.

\begin{thm} \label{thm:finiteflatfibration}
  The morphism $\psi^P:\sX^P\to\sU^{P''}$ is finite, and every irreducible
  component of $\sX^P$ surjects onto $\sU^{P''}$.
\end{thm}

\pf For dimension reasons it suffices to show that $\psi^P$ is finite.
Since $\sX^P = \sM(\cA)$ and
$\sU^P = \sM(K\angles{\sU^{P''}})$ are both affinoid, by the
rigid-analytic direct image theorem~\cite[Theorem~9.6.3/1]{bgr:nonarch} it
suffices to show that $\sX^P\to\sU^{P''}$ is proper in the sense
of~\cite[\S9.6.2]{bgr:nonarch}.  In fact we will show that
$\sX^P\Subset_{\sU^{P''}}\sX^P$, i.e.,\ that there exist affinoid
generators $f_1,\ldots,f_r$ for $\cA$ over $K\angles{\sU^{P''}}$ such that
$|f_1|_{\sup},\ldots,|f_r|_{\sup} < 1$.

Choosing bases for $N'$ and $N''$, we obtain isomorphisms
$N'_\R\cong\R^{n'}$, $N''_\R\cong\R^{d}$, and
$N_\R\cong\R^{n'}\times\R^{d}$.  
Translating by an element of $\T(K)$,
we may and do assume that $P\subset\{0\}\times N''_\R$ (so $P = P''$).
For $\epsilon\in G$ with $\epsilon > 0$ we let 
$I_\epsilon\subset N'_\R$ be the cube $[-\epsilon,\epsilon]^{n'}$, so 
$I_\epsilon$ is a integral $G$-affine polytope in $N_\R'$, and
$P_\epsilon\coloneq I_\epsilon\times P''$ is a integral $G$-affine polytope in
$N_\R = N'_\R\times N''_\R$ containing $P$.  Since $\tau$ is a maximal
face we have $P_\epsilon\cap\Trop(X) = P$ for small $\epsilon$; we fix
such an $\epsilon$ as well as an element $e\in K$ with 
$\val(e) = \epsilon$.  The polytopal subdomain
$\sU^{I_\epsilon}\subset(\T')^\an$ is a product of annuli of inner radius
$|e|$ and outer radius $|e|\inv$, so if $u_1,\ldots,u_{n'}$ is a basis for
$M'$ then $\{ex^{\pm u_1},\ldots,ex^{\pm u_{n'}}\}$ is a set of affinoid
generators for $K\angles{\sU^{I_\epsilon}}$.  Since
$\sU^{P_\epsilon} = \sU^{I_\epsilon}\times_K\sU^{P''}$, it follows that
$\{ex^{\pm u_1},\ldots,ex^{\pm u_{n'}}\}$ is a set of affinoid generators
for $K\angles{\sU^{P_\epsilon}}$ over $K\angles{\sU^{P''}}$.  Since 
$P_\epsilon\cap\Trop(X) = P$ we have $X^\an\cap\sU^{P_\epsilon}=\sX^P$, so
$\{ex^{\pm u_1},\ldots,ex^{\pm u_{n'}}\}$ can be regarded as a set of
affinoid generators for $\cA$ over $K\angles{\sU^{P''}}$.  But by
construction $|x^{u_i}(x)| = 1$ for all $x\in\sX^P$ and all
$i=1,\ldots,n'$, so $|ex^{u_i}(x)| = |e| < 1$.  This proves that
$\psi^P$ is finite.\qed

It follows from Theorem~\ref{thm:finiteflatfibration} and
Remarks~\ref{rem:pure.degree}(1) and~\ref{rem:analytic.puredegree}(1) that
$\psi^P$ has a (pure) degree.

\begin{prop} \label{prop:deg.from.mrel}
  In the situation of~\parref{par:finiteflatfibration}, let
  $\sY\subset\sX^P$ be a union of connected components and let 
  $w\in P$.   Then
  \[ [ \sY : \sU^{P''} ] = \sum_{x\in\sY\cap\trop\inv(w)} m_\rel(x). \]
\end{prop}

\pf Extending the ground field if necessary, we assume that $w\in N_G$.
Let $w''$ be the image of $w$ in $N_G''$.
Since $\sY\cap\sX^w\to\sU^{w''}$ is obtained by flat base change from 
$\sY\to\sU^{P''}$ we may replace
$P$ by $w$ and $P''$ by $w''$ to assume 
that $\sY \subset\sX^w$ (cf.\ Proposition~\ref{prop:degree.cover}(1)).  
Let $\fY$ be the canonical model of $Y$.
The canonical reduction $\bar\fY$ of $\sY$ is a union of connected
components of $\bar\fX{}^w_\can$, so for $x\in\sY$ the relative multiplicity
$m_\rel(x)$ is nonzero if and only if $\red(x)$ is the generic point of an
irreducible component $\bar\fC$ of $\bar\fY$, in which case 
$m_\rel(x) = [\bar\fC:\im(\bar\fC)]$ where $\im(\bar\fC)$ is the image
of $\bar\fC$ in $\inn_w(X)$.
Noting that $\bar\fU{}^{w''}$ is an integral domain and $\bar\fY$ is reduced,
applying the projection formula~\parref{prop:projection.formula} to 
$\fY\to\fU^{w''}$ yields
\[ [\sY:\sU^{w''}] = \sum_{\fC\subset\bar\fY} [\fC:\bar\fU{}^{w''}]. \]
Since $\bar\fU{}^{w''} = (\bar\T{}'')^{w''}$ and
$\inn_w(X)\cong\bar D\times(\bar\T{}'')^{w''}$ for some
dimension-zero scheme $\bar D\subset(\bar\T{}')^{w'}$ 
(cf.\ Remark~\ref{rem:split.torus}), the reduced space underlying any
irreducible component of $\inn_w(X)$ is isomorphic to
$(\bar\T{}'')^{w''}$.  Therefore $[\im(\bar\fC):\bar\fU{}^{w''}]=1$ for any
irreducible component $\bar\fC\subset\bar\fY$, so
$[\bar\fC:\bar\fU{}^{w''}] = [\bar\fC:\im(\bar\fC)]$ and the proposition
follows.\qed

\begin{cor}[Tropical projection formula]
  \label{cor:finiteflatfibration}
  In the situation of Theorem~\ref{thm:finiteflatfibration}, the degree of
  $\psi^P: \sX^P\to\sU^{P''}$ is equal to $m_{\Trop}(w)$ for any
  $w\in P$.
\end{cor}

\pf Assuming that $P = w$ and $P'' = w''$ as in the proof
of Proposition~\ref{prop:deg.from.mrel}, the result follows immediately
from Propositions~\ref{prop:deg.from.mrel}
and~\ref{prop:relativemults}.\qed

\begin{rem}
  The tropical projection formula is an equality of the degree of the
  morphism $\sX^P\to\sU^{P''}$ (a morphism on the generic fiber) with
  the degree of a morphism $\bar\fX\s{w}\to\bar\fU\s{w''}$ (a morphism on
  the special fiber).  It is conceptually very close to
  the projection formula as stated in
  Proposition~\ref{prop:projection.formula}, as indeed that is the main
  tool used in its proof; it is for this reason that we call it the
  tropical projection formula.
\end{rem}

\section{The tropicalization of a nonarchimedean analytic curve}
\label{section:curves}

In this section we freely use the definitions and notations
from~\cite{bpr:berk_curves}.  In particular, for $a,b\in K^\times$ we have
$\bS(a,b) = \{ t\in\bG_m^\an~:~|a|\leq|t|\leq|b|\}$, a closed annulus; we also
set $\bS(a) = \bS(a,1)$ and $\bS(1) = \{t\in\bG_m^\an~:~|t|=1\}$, the closed
annulus of modulus zero.

\paragraph[The setup] \label{par:curves.setup}
Throughout this section $X$ denotes a smooth connected algebraic curve realized as a
closed subscheme of a torus $\T$, $\hat X$ is the smooth completion of
$X$, and $D = \hat X(K)\setminus X(K)$ is the set of punctures.  We
will denote a choice of semistable vertex set for $X$ by $V$, and we let
$\Sigma = \Sigma(X, V)$ be the associated skeleton.  See~\cite[\S3]{bpr:berk_curves}.

If we choose a basis for $M$, then we obtain isomorphisms $N_\R\cong\R^n$ and
$K[\T]\cong K[x_1^{\pm 1},\ldots,x_n^{\pm 1}]$; if 
$f_i\in K[X]^\times$ is the image of $x_i$ then 
\begin{equation}
  \label{eq:trop.explicit}
  \trop(\|\cdot\|) = (-\log\|f_1\|,\ldots,-\log\|f_n\|).
\end{equation}

\paragraph[Compatible polyhedral structures]
The tropicalization of $X$ is a polyhedral complex of pure dimension $1$
in $N_\R$.  We can regard $\Trop(X)$ as a dimension-$1$ abstract
$G$-rational polyhedral complex where the metric on the edges is given by
the lattice length, i.e.\ the length in the direction of a primitive
lattice vector in $N$.  

Recall the following consequence of the 
Slope Formula~\cite[Theorem~5.15]{bpr:berk_curves}:

\begin{lem} \label{lem:integer.expansion}
  Let $e$ be an edge of $\Sigma$ and let $f\in K[X]^\times$.  The
  map $\|\cdot\|\mapsto-\log\|f\|:  X^\berk\to\R$ restricts to a $G$-affine
  linear function from $e$ to $\R$ with integer slope.  
\end{lem}

Since the $G$-rational points of an edge $e$ of $\Sigma$ are exactly the
type-$2$ points of $X^\an$ contained in $e$
(cf. \cite[(3.12)]{bpr:berk_curves}), it
follows from~\eqref{eq:trop.explicit} and the
$G$-rationality of $-\log\|f\|$ as above that $\trop$ maps type-$2$ points
into $N_G$.

\begin{prop} \label{prop:compatible.polyhedral}
  \begin{enumerate}
  \item The map $\trop: X^\berk\surject\Trop(X)$ factors through the
    retraction $\tau_\Sigma: X^\berk\surject\Sigma$.
  \item We can choose $V$ and a polyhedral complex structure on $\Trop(X)$
    as in Theorem~\ref{thm:polyhedral.decomp}
    such that $\trop: \Sigma\to \Trop(X)$ is a morphism of dimension-$1$ abstract
    $G$-rational polyhedral complexes.
  \end{enumerate}
\end{prop}

\pf The first part follows from~\eqref{eq:trop.explicit} and the
Slope Formula~\cite[Theorem~5.15]{bpr:berk_curves}.
Let $e$ be an edge of $\Sigma$.
It follows from Lemma~\ref{lem:integer.expansion} as applied to
$f_1,\ldots,f_n$ that $\trop$ restricts to an expansion by an 
integer multiple with respect to the intrinsic metric on $e$ and the
lattice length on its image.
Hence there exist refinements of the polyhedral structures
on $\Sigma$ and on $\Trop(X)$ such that $\trop: \Sigma\to\Trop(X)$
becomes a morphism of dimension-$1$ abstract $G$-rational polyhedral
complexes.  By \cite[Proposition 3.13(2)]{bpr:berk_curves}, any
refinement of $\Sigma$ is also a skeleton of $X$.\qed

From now on we assume that our skeleton $\Sigma$ of $ X$ and our choice
of polyhedral structure on $\Trop(X)$ are compatible in the sense of
Proposition~\ref{prop:compatible.polyhedral}(2).  If $e\subset\Sigma$
is an interval contained in an edge we let $\ell_{\Berk}(e)$ be its length
with respect to the skeletal metric, and if $e'\subset\Trop(X)$ is an
interval contained in an edge we let $\ell_{\Trop}(e')$ be its lattice
length.  As a consequence of Proposition~\ref{prop:compatible.polyhedral}(2), if
$\trop(e) = e'$ then $\ell_{\Trop}(e')$ is an integer multiple of
$\ell_{\Berk}(e)$. 

\begin{defn}
  Let $e\subset\Sigma$ be an edge and let $e'\subset\Trop(X)$ be its
  image.  We define the \emph{expansion factor of $e$} to be the unique integer
  $m_\rel(e)\in\Z_{\geq 0}$ such that
  \[ \ell_{\Trop}(\trop(\tilde{e})) = m_\rel(e)\cdot\ell_{\berk}(\tilde{e}) \]
  for any finite-length segment $\tilde{e}$ contained in $e$.
\end{defn}

\begin{rem} \label{rem:mult.is.gcd}
  For $u\in M$ let $f_u\in K[X]^\times$ be the image of the character
  $x^u\in K[M]$.  Let $e$ be an edge of $\Sigma$ and let
  $s_u\in\Z_{\geq 0}$ be the absolute value of the slope of $-\log\|f_u\|$
  on $e$.  It follows easily from the definitions that
  $m_\rel(e) = \gcd\{s_u~:~u\in M\}$.  More concretely, let
  $u_1,\ldots,u_n$ be a basis for $M$ and let $f_i\in K[X]^\times$ be
  the image of $x^{u_i}$, so
  $\trop(\|\cdot\|) = (-\log\|f_1\|,\ldots,-\log\|f_n\|)$.
  Let $s_i\in\Z_{\geq 0}$ be the absolute value of the
  slope of $-\log\|f_i\|$ on $e$.  Then
  \begin{equation} \label{eq:mrel.is.gcd.slopes}
    m_\rel(e) = \gcd(s_1,\ldots,s_n).
  \end{equation}
\end{rem}

\paragraph \label{par:compat.mult.setup} We now come to one of the key
results of this section.  Let $e\subset\Sigma$ be a bounded edge,%
\footnote{If $e$ is an infinite ray, we can compute $m_{\rel}(e)$ by
  refining the polyhedral structure on $\Sigma$ so that $m_{\rel}(e)$
  coincides with $m_{\rel}(\tilde{e})$ for some bounded edge $\tilde{e}$
  in the refinement.  Thus we can assume without loss of generality in
  Theorem~\ref{thm:compatible.multiplicities} that $e$ is bounded.}
and assume that
$e' = \trop(e)$ is an edge of $\Trop(X)$ (as opposed to a vertex).  
The inverse image of the interior of $e'$ under $\trop$ is a disjoint
union of open annuli, one of whose skeleta is the interior of $e$.  Hence
if $x\in e$ is the unique
point mapping to some $w \in \relint(e')\cap N_G$ then 
$\sY_x \coloneq \tau_\Sigma\inv(x)\cong\bS(1)$ is a 
connected component of $\sX^w = \sU^w\cap X^\an$ 
(Definition~\ref{defn:polyhedral.models}).  Let
$W$ be the affine span of $e'$, let $\T'$ be the torus transverse to
$W$~\parref{par:polyhedral.structures}, let $\T''$ be the torus parallel
to $W$ (Remark~\ref{rem:split.torus}), and choose a splitting
$\T\to\T''$.
We have a finite surjective morphism  
$\sX^w\to\sU^{w''}$ by Theorem~\ref{thm:finiteflatfibration},
where $w''$ is the image of $w$ in $N_G''$.  

For $y\in e$ the relative multiplicity $m_\rel(y)$ was defined
in~\parref{defn:m.rel.xi}.  

\begin{thm}[Compatibility of multiplicities]
  \label{thm:compatible.multiplicities}
  With the above notation, 
  \[ m_\rel(e) = [ \sY_x : \sU^{w''} ] = m_\rel(x). \]
  Moreover $m_\rel(e) = m_\rel(y)$ for any $y$ in the interior of $e$
  (even if $\trop(y)\notin N_G$).
\end{thm}

\pf Let $P$ be a $G$-rational closed interval contained in the interior of
$e'$ and containing $w$.  Then 
\[ \sY\coloneq \tau_\Sigma\inv(\trop\inv(P)\cap e)\cong\bS(a,b) \]
is a connected component of $\sX^P$; it is a closed annulus of nonzero
modulus with skeleton $\trop\inv(P)\cap e$.  Let $P''$ be the image of $P$
in $N_G''$, so $\ell_{\Trop}(P'') = \ell_{\Trop}(P)$.
By Theorem~\ref{thm:finiteflatfibration} the morphism  
$\sY\to\sU^{P''}$ is finite and surjective, and
$[\sY:\sU^{P''}]=[\sY_x:\sU^{w''}]$ 
by Proposition~\ref{prop:degree.cover}(1).  The 
torus $\T''$ is one-dimensional, so $\sU^{P''} = \trop\inv(P'')$ is an
annulus and $P''$ is by definition the skeleton of $\sU^{P''}$.  We have
\[ m_\rel(e)\cdot\ell_{\berk}(\trop\inv(P)\cap e) = \ell_{\Trop}(P) =
\ell_{\Trop}(P'') = [\sY:\sU^{P''}]\cdot \ell_{\berk}(\trop\inv(P)\cap e), \]
where the final equality is by \cite[Corollary~2.6]{bpr:berk_curves}.
Since $x$ is the unique Shilov boundary point of $\sY_x\cong\bS(1)$,
the equality $[\sY_x:\sU^{w''}] = m_\rel(x)$ is a consequence of
Proposition~\ref{prop:deg.from.mrel}.

By a standard argument involving extension of the ground field
(cf.\ Lemma~\ref{lem:multiplicities.extension.scalars}
and~\parref{par:base.change}), the second statement follows from the first.\qed

\begin{cor}\label{cor:Berkovichmultiplicities}
Fix an edge $e'$ of $\Trop(X)$ and let $e_1,\ldots,e_r$ be the edges of
$\Sigma$ mapping homeomorphically onto $e'$.  Then
\[m_{\Trop}(e') = \sum_{i=1}^r m_{\rel}(e_i). \]
\end{cor}

\pf This follows immediately from
Theorem~\ref{thm:compatible.multiplicities} 
and Proposition~\ref{prop:relativemults}.\qed

\begin{rem}
  With the notation in Corollary~\ref{cor:Berkovichmultiplicities}, 
  let $w$ be a $G$-rational point contained in the relative interior of
  $e'$.   The affinoid space $\trop\inv(w)$ is isomorphic to
  $\Djunion_{i=1}^r \bS(1)$ as in~\parref{par:compat.mult.setup}.
  Since the canonical reduction of $\bS(1)$ is isomorphic to $\G_m$,
  the integer $r$ is equal to the number of
  irreducible components in the canonical reduction of $\sX^w$.
\end{rem}

\begin{cor}
\label{cor:multiplicityoneisometry}
If $e'$ is an edge of $\Trop(X)$ and $m_{\Trop}(e')=1$, then there is a
unique edge $e$ in $\Sigma$ such that $\trop$ maps $e$ homeomorphically
and isometrically onto $e'$.  The edge $e$ is in fact the unique geodesic
segment (or ray) in $X^\Berk$ which is mapped homeomorphically by $\trop$ onto
$e'$. 
\end{cor}

\begin{cor} \label{cor:skeleton.multiplicities}
  Let $x\in\HH(X^\an)$.  Then $m_\rel(x) > 0$ if and only if $x$ belongs to an
  edge of $\Sigma$ mapping homeomorphically onto its image via $\trop$.
\end{cor}

\pf Suppose that $x$ is contained in an edge $e\subset\Sigma$ mapping
homeomorphically onto its image $e' = \trop(e)$.  If $x$ is in the
interior of $e$ then $m_\rel(x) = m_\rel(e) > 0$ 
by Theorem~\ref{thm:compatible.multiplicities}.  Otherwise 
$w = \trop(x)\in N_G$, and $x$ is contained in the limit boundary 
$\del_{\lim}\sX^w$ of $\sX^w = \trop\inv(w)$ because it is a limit of points
of $e$ which are not contained in $\sX^w$.  Since
$\del_{\lim}\sX^w$ is the Shilov boundary of $\sX^w$, by definition we have
$m_\rel(x) > 0$.  

Now suppose that $x$ is not contained in an edge of $\Sigma$ mapping
homeomorphically onto its image.  If $x$ has type $4$ then $m_\rel(x)=0$
by definition.  Otherwise by \cite[Corollary 5.1]{bpr:berk_curves}
we can enlarge $\Sigma$ if necessary to assume that $x\in\Sigma$ .  
Recall that every edge of $\Sigma$ maps homeomorphically onto its image or
is crushed to a vertex of $\Trop(X)$.  By hypothesis all edges containing
$x$ are crushed to a vertex of $\Trop(X)$, so $w = \trop(x)\in N_G$.
Hence there is an open neighborhood $U$ of $x$ in $\Sigma$ contained in 
$\sX^w$.  Then $\tau_\Sigma\inv(U)$ is a neighborhood of $x$ in $X^\an$
contained in $\sX^w$, so $x\notin\del\sX^w$ and hence $m_\rel(x)=0$.\qed

\paragraphtoc[Slopes as orders of vanishing]
\label{par:slopes.orders.vanishing}
The Slope Formula~\cite[Theorem~5.15]{bpr:berk_curves}
provides a useful
interpretation of the quantities $s_i$ appearing in
Remark~\ref{rem:mult.is.gcd} in terms of orders of vanishing.  Assume that
$V$ is a strongly semistable vertex set~\cite[\S3]{bpr:berk_curves} of $\hat X$ (in addition to being
a semistable vertex set of $X$).  Let $\fX$ be the strongly semistable
formal model of $\hat X$ associated to $V$,
let $x\in V$, and let $\bar\fC$ be
the irreducible component of $\bar\fX$ whose generic point is $\red(x)$.
Let $e$ be an edge of $\Sigma$ adjacent to $x$ and let $\xi\in\bar\fC(k)$ be
the reduction of the interior of $e$.  The Slope Formula~\cite[Theorem~5.15]{bpr:berk_curves}
says that if
$f$ is a nonzero rational function on $\hat X$ then the 
slope $s$ of $-\log|f|$ along $e$ (in the direction away from $x$) 
is equal to $\ord_\xi(\td f_x)$.  One
can use this fact to give a simple proof of the well-known balancing
formula for tropical curves:

\begin{thm}[The balancing formula for tropical curves]
\label{thm:balancingformula}
Let $w$ be a vertex of $\Trop(X)$ and let $\vec{v}_1,\ldots,\vec{v}_t$ be
the primitive integer tangent directions at $w$ corresponding to the
various edges $e'_1,\ldots,e'_t$ incident to $w$.  Then 
$\sum_{j=1}^t m_{\Trop}(e'_j) \vec{v}_j = 0$.
\end{thm}

\pf We use the setup in~\parref{par:slopes.orders.vanishing}.
Let $f_1,\ldots,f_n\in K[X]^\times$ be the coordinate functions as
in~\parref{rem:mult.is.gcd} and let $F_i = -\log|f_i|$.
Let $x \in \Sigma \cap \trop^{-1}(w)$ be a vertex.
By the Slope Formula~\cite[Theorem~5.15]{bpr:berk_curves} we have 
$0 = \sum_{v\in T_x} d_v F_i(x)$ for each $i=1,\ldots,n$.
Since we are assuming that $V$ is a strongly semistable vertex set, each
tangent direction $v\in T_x$ along which some $F_i$ has nonzero slope
is represented by a unique edge 
$e_v = [x,y_v]$ of $\Sigma$ adjoining $x$, and $d_v F_i(x)$ is just the
slope of $-\log|f_i|$ along $e_v$.
If $\trop(e_v) = \{w\}$ then $d_v F_i(x) = 0$ for all
$i$, and otherwise 
\[ d_v F_i(x) = \frac{\log|f_i(x)| - \log|f_i(y_v)|}{\ell_\an(e_v)} 
= m_\rel(e_v)\, \frac{\log|f_i(x)| -
  \log|f_i(y_v)|}{\ell_{\Trop}(\trop(e_v))}. \] 
By Corollary~\ref{cor:Berkovichmultiplicities}, for each $i$ we have 
\[
\begin{aligned}
\sum_{j=1}^t m_{\Trop}(e_j') (\vec{v}_j)_i &=  
\sum_{j=1}^t\bigg( \sum_{e = [x,y] \isom e_j'} m_{\rel}(e)\bigg) \,
\frac{\log|f_i(x)| - \log|f_i(y)|}{\ell_{\Trop}(e_j')}   \\
&= \sum_{x \mapsto w} \sum_{v\in T_x} d_v F_i(x) = 0, \\
\end{aligned}
\]
which implies the result.\qed

\paragraphtoc[Faithful representations]
If $Y_\Delta$ is a proper toric variety with dense torus $\T$ then
$X\inject\T$ extends in a unique way to a morphism 
$\iota: \hat X\to Y_\Delta$, which is a closed immersion for suitable
$Y_\Delta$.  The 
intersection $\hat X^\an\cap(Y_\Delta^\an\setminus\T^\an)$ is the finite
set of type-$1$ points $D = \hat X^\an\setminus X^\an$.  We write
$\trop(\iota): \hat X^\an\to N_\R(\Delta)$ for the induced tropicalization
map, and we set
$\Trop(\hat X, \iota) = \trop(\iota)(\hat X^\an)\subset N_\R(\Delta)$.

\subparagraph \label{par:order.embeddings}
Let $Y_\Delta,Y_{\Delta'}$ be toric varieties with dense tori $\T,\T'$ and
let $\iota:\hat X\inject Y_\Delta$ and $\iota':\hat X\inject Y_{\Delta'}$
be closed immersions whose images meet the dense torus.  We say that
$\iota'$ \emph{dominates} $\iota$ and we write $\iota'\geq\iota$ provided
that there exists a morphism $\psi:Y_{\Delta'}\to Y_\Delta$ of toric
varieties such that $\psi\circ\iota'=\iota$.  In this case we have an
induced morphism 
$\Trop(\psi):\Trop(\hat X,\iota')\to\Trop(\hat X,\iota)$ making the
triangle
\[\xymatrix @=.2in{
  & {\hat X^\an} \ar[dl]_{\trop(\iota')} \ar[dr]^{\trop(\iota)} & \\
  {\Trop(\hat X,\iota')} \ar[rr]_{\Trop(\psi)} & & 
  {\Trop(\hat X,\iota)}
}\]
commute.  Since $\trop(\iota)$ and $\trop(\iota')$ are surjective, the map
$\Trop(\psi)$ is independent of the choice of $\psi$, so the set of
`tropicalizations of toric embeddings' is a filtered inverse system.

\subparagraph
By a \emph{finite subgraph} of $\hat X^\an$ we mean a \emph{connected} compact
subgraph of a skeleton of $\hat X$.  Any finite union of geodesic segments
in $\HH_\circ(\Xhat^{\an})$ is contained in a skeleton by
\cite[Corollary 5.10]{bpr:berk_curves},
so we can equivalently define a
finite subgraph of $\hat X^\an$ to be an isometric embedding of a finite
connected metric graph $\Gamma$ into $\HH_\circ(\Xhat^{\an})$.
Let $\hat X\inject Y_\Delta$ be a closed immersion into a toric variety
with dense torus $\T$ such that $\hat X\cap\T\neq\emptyset$.
We say that a finite subgraph $\Gamma$ of $\hat X^\an$
is {\em faithfully represented} by $\trop:\hat X^\an\to N_\RR(\Delta)$ if
$\trop$ maps $\Gamma$ homeomorphically and isometrically onto its image
$\Gamma'$ (which is contained in $N_\R$).  
We say that $\trop$ is {\em faithful} if it faithfully represents a
skeleton $\Sigma$ of $\hat X$.

\begin{rem} \label{rem:local.isometry}
  When considering a closed connected subset $\Gamma$ of $\HH(X^\an)$ or
  $\Gamma'$ of $\Trop(X)$, we will always implicitly endow it with the
  shortest-path metric.  In general this is \emph{not} the same as the
  metric on $\Gamma$ (resp.\ $\Gamma'$) induced by the (shortest-path)
  metric on $\HH(X^\an)$ (resp.\ $\Trop(X)$).  With this convention,
  $\Gamma$ (resp.\ $\Gamma'$) is a \emph{length space} in the sense
  of~\cite[Definition~2.1.2]{papadopoulos:metric_spaces}, so any
  homeomorphism $\Gamma\to\Gamma'$ which is a local isometry is
  automatically an isometry by Corollary~3.4.6 of loc.\ cit.
  This will be used several times in what follows.
\end{rem}

The following result shows that if $\Gamma$ is faithfully represented by a
given tropicalization, then it is also faithfully represented by all
`larger' tropicalizations.   

\begin{lem} \label{lem:faithfulnessstablilizes}
Let $\iota : \Xhat \into Y_\Delta$ and $\iota' : \Xhat \into Y_{\Delta'}$ be closed
immersions of $\Xhat$ into toric varieties whose images meet the dense
torus and such that $\iota'\geq\iota$. 
If a finite subgraph $\Gamma$ of $\Xhat^{\an}$ is faithfully represented by 
$\trop(\iota):\hat X^\an\to N_\RR(\Delta)$ 
then $\Gamma$ is faithfully represented by 
$\trop(\iota'):\hat X^\an\to N_\RR(\Delta')$. 
\end{lem}

\pf Without loss of generality, we may replace $\Xhat$ by $X = \Xhat
\smallsetminus D$ and assume that $Y_{\Delta'}$ and $Y_\Delta$ are tori 
with $\iota = (f_1,\ldots,f_n)$ and $\iota' = (f_1,\ldots,f_m)$ for some 
$m \geq n$.
The result is now clear from~\eqref{eq:mrel.is.gcd.slopes}.\qed

We will show in Theorem~\ref{thm:exists.faithful.trop} below that any
finite subgraph of $\hat X^\an$ is faithfully represented by some
tropicalization.  First we need two lemmas.

\begin{lem}
  \label{lem:multone}
  Let $e$ be an edge of a skeleton $\Sigma$ of $\hat X$ with distinct
  endpoints $x,y$.  There exists a nonzero meromorphic function $f$ on
  $\hat X$ such that $F = -\log |f|$ has the following properties:
  \begin{enumerate}
  \item $F\geq 0$ on $e$, and $F(x) = F(y) = 0$.
  \item There exist (not necessarily distinct)
    type-$2$ points $x',y'$ in the interior $e^\circ$ of
    $e$ such that $\rho(x,x') = \rho(y,y')$, such that $F$ has slope $\pm 1$ on 
    $[x,x']$ and $[y,y']$, and such that $F$ is constant on $[x',y']$, 
    as shown in Figure~\ref{fig:graphofF}.
  \end{enumerate}
\end{lem}

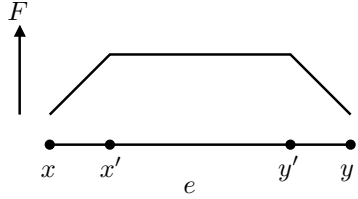
\genericfig[ht]{graphofF}{The graph of the function 
  $F = -\log |f|:e\to\R_{\geq 0}$
  constructed in Lemma~\ref{lem:multone}.}

\pf We may assume without loss of generality that
$\Sigma = \Sigma(\hat X, V(\fX))$ for a strongly
semistable formal model $\fX$ of $\hat X$~\cite[\S4]{bpr:berk_curves}.
For each irreducible component $\bar\fC_\nu$ of $\bar\fX$, a simple
argument using the Riemann-Roch theorem allows us to choose a
rational function $\td{f}_\nu$ on $\bar\fC_\nu$ which vanishes to order
$1$ at every singular point of $\bar\fX$ lying on $\bar\fC_{\nu}$.
By~\cite[Corollary~3.8]{bosch_lutk:uniformization} 
there exists a nonzero rational function 
$f$ on $\hat X$ whose poles all reduce to smooth points of $\bar\fX$
and which induces the rational function $\td{f}_{\nu}$ on each
irreducible component $\bar\fC_{\nu}$ of $\bar\fX$ 
(the gluing condition from loc.\ cit.\ is trivially   in this
situation).  The function $f$ constructed in the proof of 
loc.\ cit. is defined on an affinoid domain $U$ of $\hat X^\an$
containing $x$ and $y$, and $\td f_x,\td f_y$ are the
restrictions of the residue of $f$ in the canonical reduction of $U$.
Therefore we have $|f(x)| = |f(y)| = 1$.  Since $\{x,y\}$ is
the Shilov boundary of $\tau_\Sigma\inv(e)$, this proves~(1).

By~\cite[Theorem~5.15(3)]{bpr:berk_curves}
the outgoing slope of $F$ at $x$ or $y$ in the 
direction of $e$ is $1$.  Let 
$\xi$ be the singular point of $\bar\fX$ whose formal fiber is 
$\tau_\Sigma\inv(e^\circ)$.  Since $f$ has no poles on the formal fibers above
singular points, $f$ restricts to an analytic function on the open annulus 
$\tau_\Sigma\inv(e^\circ)$.  Part~(2) now follows from~(1) and
\cite[Proposition 2.10(1)]{bpr:berk_curves}.\qed

\begin{lem} \label{lem:multone.annulus}
  Let $A\subset\hat X^\an$ be an affinoid domain isomorphic to a closed
  annulus $\bS(a)$ with nonzero modulus.  There exists a nonzero
  meromorphic function $f$ on $\hat X$ such that $F = -\log |f|$ is linear
  with slope $\pm1$ on $\Sigma(A)$. 
\end{lem}

\pf Choose an identification of $A$ with 
$\bS(a) = \sM(K\angles{at\inv,t})$.  
By~\cite[Th\'eor\`eme~4,~\S2.4]{fresnel_matignon:dim_1}
the ring of meromorphic functions on $\hat X$ which are regular on $A$ is
dense in $K\angles{at\inv,t}$.  Hence there exists a meromorphic function 
$f$ on $\hat X$ such that $f\in K\angles{at\inv,t}$ and 
$|f-t|_{\sup} < 1$.  It follows from 
\cite[Proposition 2.2]{bpr:berk_curves}
that $f$ is also a parameter for
the annulus $A$, so $-\log |f|$ is linear with slope $\pm 1$ on
$\Sigma(A)$.\qed

\begin{thm} \label{thm:exists.faithful.trop}
  If $\Gamma$ is any finite subgraph of $\Xhat^{\an}$ then there is a closed
  immersion $\hat X\inject Y_\Delta$ of $\hat X$ into a quasiprojective toric variety $Y_\Delta$ such
  that $\trop:\hat X^\an\to N_\RR(\Delta)$ faithfully represents $\Gamma$.
  In particular, there exists a faithful tropicalization.
\end{thm}

\pf Since $\Gamma$ is by definition contained in a skeleton $\Sigma$, we
may assume that $\Gamma = \Sigma$.  Taking a refinement of $\Sigma$ if
necessary, we assume without loss of generality that $\Sigma$ does not
have any loop edges.  We claim that after possibly refining $\Sigma$ further,
for each edge $e\subset\Sigma$ there exists a nonzero
meromorphic function $f$ on $\hat X$ such that $\log|f|$ has slope $\pm 1$
on $e$.

Let $e = [x,y]$ be an edge of $\Sigma$, and let $f$ and $x',y'\in e^\circ$ be
as in Lemma~\ref{lem:multone}.  Then $[x,x']$ and $[y',y]$ are edges in a
refinement of $\Sigma$, and $\log|f|$ has slope $\pm 1$ on 
$[x,x']$ and $[y',y]$.  If $x'=y'$ then we are done with $e$; otherwise we
let $e'= [x',y']$.  By construction $e'\subset e^\circ$, so 
$\tau_\Sigma\inv(e')$ is a closed annulus of nonzero modulus, and we may apply
Lemma~\ref{lem:multone.annulus} to find $f'$ such that $\log|f'|$ has
slope $\pm 1$ on $e'$.  This proves the claim.

By~\eqref{eq:mrel.is.gcd.slopes}, if $\Phi = \{f_1,\ldots,f_r\}$ is any collection of
meromorphic functions on $X$ such that (a) for each edge $e$ of $\Sigma$ there
is an $i$ such that $\log|f_i|$ has slope $\pm 1$ on $e$, and (b)
$\phi = (f_1,\ldots,f_r)$ induces a closed immersion of a dense open
subscheme $X$ of $\hat X$ into a torus $\T\cong\G_m^r$, then 
$\trop\circ\phi$ maps each edge of $\Sigma$ isometrically onto its
image.  Since $\phi$ extends to a closed
immersion $\hat X\inject Y_\Delta$ into a 
suitable compactification $Y_\Delta$ of $\T$, it only remains to show that we can
enlarge $\Phi$ so that $\phi|_\Sigma$ is \emph{injective}; then
$\trop\circ\phi$ maps $\Sigma$ isometrically onto its image
by Remark~\ref{rem:local.isometry}.  

Let $e$ be an edge of $\Sigma$.  Since $\Sigma$ has at least two
edges, 
the proof of \cite[Theorem 4.11]{bpr:berk_curves} (specifically Case 1 in (4.15.1) of {\em loc. cit.})
shows that
$\tau_\Sigma\inv(e)$ is an affinoid domain in $\hat X^\an$.
By~\cite[Th\'eor\`eme~1,~\S1.4]{fresnel_matignon:dim_1}, there is a
meromorphic function $f$ on $\hat X$ such that 
$\tau_\Sigma\inv(e) = \{x\in\hat X^\an~:~|f(x)|\leq 1\}$.  Adding such an $f$ to
$\Phi$ for every edge $e$, we may assume that $\trop\circ\phi$ is
injective on $\Sigma\setminus V$, i.e.,\ that if 
$\trop(\phi(x)) = \trop(\phi(y))$ for $x,y\in\Sigma$ then $x,y$ are
vertices.  By the definition of $\hat X^\an$, if $x,y\in\hat X^\an$ are
distinct points then there exists a meromorphic function $f$ on $\hat
X^\an$ such that $|f(x)|\neq|f(y)|$.  Adding such $f$ to $\Phi$ for every pair
of vertices yields a faithful tropicalization.\qed

We obtain the following theorem as a consequence:

\begin{thm} \label{thm:MainThm4}
  Let $\Gamma$ be a finite subgraph of $\hat X$.
  Then there exists a quasiprojective toric embedding 
  $\iota : \hat X \into Y_{\Delta}$ such that for every 
  quasiprojective toric embedding $\iota' : \hat X \into Y_{\Delta'}$ 
  with $\iota' \geq \iota$, the tropicalization map 
  $\trop(\iota') : \hat X^{\an} \to N_{\RR}(\Delta')$ 
  maps $\Gamma$ homeomorphically and isometrically onto its image.
\end{thm}

\pf By Theorem~\ref{thm:exists.faithful.trop}, there exists a closed
embedding $\iota$ such that $\trop(\iota)$ maps $\Gamma$ homeomorphically and
isometrically onto its image.  By Lemma~\ref{lem:faithfulnessstablilizes},
the same property holds for any closed embedding $\iota' \geq \iota$. \qed

As mentioned in the introduction, Theorem~\ref{thm:MainThm4} can be
interpreted colloquially as saying that the homeomorphism in
Theorem~\ref{thm:paynehomeo} is an isometry.

\smallskip

With a little more work, we obtain the following strengthening of Theorem~\ref{thm:MainThm4} in which
the finite metric graph $\Gamma$ is replaced by an arbitrary skeleton of $X$ (which is 
no longer required to be compact or of finite length).

\begin{thm} \label{thm:MainThm4bis}
  Let $\Sigma$ be any skeleton of $X$.
  Then there exists a quasiprojective toric embedding 
  $\iota : X \into Y_{\Delta}$ such that for every 
  quasiprojective toric embedding $\iota' : X \into Y_{\Delta'}$ 
  with $\iota' \geq \iota$, the tropicalization map 
  $\trop(\iota') : \hat X^{\an} \to N_{\RR}(\Delta')$ 
  maps $\Sigma$ homeomorphically and isometrically onto its image.
\end{thm}

\pf
Using Lemma~\ref{lem:faithfulnessstablilizes}, it suffices to prove that
there exists a closed embedding $\iota$ such that $\trop(\iota)$ maps $\Sigma$ homeomorphically and
isometrically onto its image.  

For each point $p \in D = \Xhat \smallsetminus X$, choose a pair of
relatively prime integers $m_1(p),m_2(p)$ bigger than $2g$, where $g$ is
the genus of $\Xhat$.  By the Riemann-Roch theorem, there are rational
functions $f_1^{(p)}$ and $f_2^{(p)}$ on $\Xhat$ such that $f_i^{(p)}$ has
a pole of exact order $m_i(p)$ at $p$ and no other poles for $i=1,2$.  Let
$U_p$ be an (analytic) open neighborhood of $p$ on which $f_1^{(p)}$ and $f_2^{(p)}$
have no zeros and let $U$ be the union of $U_p$ for all $p \in D$.

Let $\Gamma = \Sigma \smallsetminus (\Sigma \cap U)$.  Then $\Gamma$ is a
finite subgraph of $X$ so by Theorem~\ref{thm:exists.faithful.trop}, there
exists a closed embedding $\iota_0$ such that $\trop(\iota_0)$ maps
$\Gamma$ homeomorphically and isometrically onto its image.  We can choose
the $U_p$ such that the
complement $\Sigma \smallsetminus \Gamma$ consists of finitely many
open infinite rays $r_p$, one for each point $p \in D$.
By the Slope Formula, the absolute value of the slope of $\log|f_i^{(p)}|$
along $r_p$ is $m_i(p)$.  Since $\gcd(m_1(p),m_2(p))=1$ for all $p \in D$,
if we enhance the embedding $\iota_0$ to a larger embedding $\iota$ by
adding the coordinate functions $f_1^{(p)}$ and $f_2^{(p)}$ for all 
$p \in D$, $\trop(\iota)$ has multiplicity one along each ray $r_p$ by
Remark~\ref{rem:mult.is.gcd}.  By Lemma~\ref{lem:faithfulnessstablilizes},
$\trop(\iota)$ also has multiplicity one at every edge of $\Gamma$.  It
follows easily (as in the proof of Theorem~\ref{thm:exists.faithful.trop})
that $\trop(\iota)$ maps $\Sigma$ homeomorphically and isometrically onto
its image as desired.  \qed

\paragraphtoc[Certifying faithfulness]
It is useful to be able to certify that a given tropicalization map is faithful using only `tropical' computations.

\begin{thm}
\label{thm:SectionThm}
Let $\Gamma'$ be a compact connected subset of $\Trop(X)$
and suppose that $m_{\Trop}(w) = 1$ for all $w \in \Gamma'\cap N_G$.
Then there is a unique closed subset $\Gamma \subset \HH_\circ(X^\an)$
mapping homeomorphically onto $\Gamma'$, and this homeomorphism is an isometry.
\end{thm}

\pf Since $m_{\Trop}$ is constant along the interior of each edge of
$\Trop(X)$ by Theorem~\ref{thm:polyhedral.decomp}(1) 
(for points not contained in $N_G$ this is proved by
a standard ground field extension argument) and $\Gamma'$ is a
finite union of closed intervals, we have 
$m_{\Trop}(w) = 1$ for all $w\in\Gamma'$.  By
Proposition~\ref{prop:relativemults}, for each $w \in \Gamma'$ there is a
unique point $x = x_w \in\HH_\circ(X^\an)$ such that $\trop(x)=w$ and
$m_{\rel}(x)>0$.  Let $\Gamma = \{ x_w \; : \; w \in \Gamma' \}$.  The
natural continuous map $\trop : \Gamma \to \Gamma'$ is bijective.
It follows from Corollaries~\ref{cor:multiplicityoneisometry}
and~\ref{cor:skeleton.multiplicities} that $\Gamma$ is also a finite union of
closed intervals, hence compact.  Thus $\trop : \Gamma \to \Gamma'$, being
a continuous bijection between compact Hausdorff spaces, is a
homeomorphism.  By Corollary~\ref{cor:multiplicityoneisometry} this
homeomorphism is an isometry. 

As for the uniqueness of $\Gamma$, let $\td\Gamma$ be any
closed subset of $\HH_\circ(X^\an)$ mapping homeomorphically onto
$\Gamma'$.  Fix $w \in \Gamma'$, and let $x$ be the point in $\td\Gamma$
with $\trop(x)=w$.  Since $x$ belongs to a closed segment of $X^{\an}$
mapping homeomorphically onto its image via $\trop$ (namely the inverse
image in $\td\Gamma$ of an edge in $\Trop(X)$ containing $w$), it follows from
Corollary~\ref{cor:skeleton.multiplicities} that $m_{\rel}(x) > 0$.  Hence
$x = x_w$, so $\td\Gamma = \Gamma$.\qed

In order to apply Theorem~\ref{thm:SectionThm}, it is useful to know that
one can sometimes determine the multiplicity at a point $w \in \Trop(X)
\cap N_G$ just from the local structure of $\Trop(X)$ at $w$, i.e.,\ from
the combinatorics of $\Star(w)$.  Recall that if $\vec v_1, \ldots, \vec
v_r$ are the primitive generators of the edge directions in $\Trop(X)$ at
$w$, and $a_i$ is the tropical multiplicity of the edge corresponding to
$\vec v_i$, then the balancing condition says that $a_1 \vec v_1 + \cdots
+ a_r \vec v_r = 0$.  Now, if $Z_1, \ldots, Z_s$ are the irreducible
components of $\inn_w(X)$ then the tropicalization of each $Z_j$ (as a
subscheme of the torus torsor $\bar\T\s w$ over the trivially-valued field
$k$) is a union of rays spanned by a subset of $\{\vec v_1, \ldots, \vec
v_r \}$.  If $b_{ij}$ is the multiplicity of the ray spanned by $\vec v_i$
in $\Trop(Z_j)$ and $m_i$ is the multiplicity of $Z_i$ in $\inn_w(X)$ then
then the balancing condition implies that $b_{1j} \vec v_1
+ \cdots + b_{rj}\vec v_r = 0$, and we also have 
$m_1 b_{i1} + \cdots + m_s b_{is} = a_i$, since 
$\Trop( \inn_w(X)) = \Star_w(\Trop(X))$ by
\cite[Proposition 10.1]{speyer:uniformizing}.

\begin{thm}
\label{thm:combmult1theorem}
Let $w \in \Trop(X) \cap N_G$.  If $\Trop(X)$ is trivalent at $w$ and one
of the edges adjacent to $w$ has multiplicity one, then $w$ has multiplicity
one. 
\end{thm}

\pf
Let $\vec v_1, \vec v_2, \vec v_3$ be the primitive generators of the edge directions in $\Trop(X)$ at $w$, 
and let $a_i$ be the multiplicity of the edge in direction $\vec v_i$.  The linear span $\angles{ \vec v_1, \vec v_2, \vec v_3 }$ 
is two dimensional, since the $\vec v_i$ are distinct and satisfy the balancing condition, so any relation among them is 
a scalar multiple of the relation $a_1 \vec v_1 + a_2\vec v_2 + a_3 \vec v_3 = 0$.

Let $Z$ be an irreducible component of $\inn_w(X)$, and let $b_i$ be the multiplicity of the ray spanned by $\vec v_i$ in $\Trop(Z)$.  
Then $b_i$ is a nonnegative integer bounded above by $a_i$ and $b_1 \vec v_1 + b_2 \vec v_2 + b_3 \vec v_3 = 0$.  
This relation must be a scalar multiple of the relations given by the $a_i$, so there is a positive rational number $\lambda \leq 1$ 
such that $b_i = \lambda a_i$ for all $i$.  
If some $a_i$ is one, then $\lambda$ must also be one and $b_i = a_i$ for all $i$.  Since $a_i$ is the sum of the multiplicities 
of the ray spanned by $\vec v_i$ in the tropicalizations of the components of $\inn_w(X)$, it follows that $\inn_w(X)$ has 
no other components, and $w$ has multiplicity one.\qed

\begin{rem}
Initial degenerations at interior points of an edge of multiplicity $1$ are always smooth, since they are isomorphic to $\Gm$ by
Theorem~\ref{thm:polyhedral.decomp}.
\end{rem}

\begin{rem}
\label{rem:combmult1remark}
There are other natural combinatorial conditions which can guarantee multiplicity one at a point $w \in N_G$ of $\Trop(X)$
or smoothness of the corresponding initial degeneration. 
In the case of curves, for example, the argument above works more generally if $\Trop(X)$ is $r$-valent, 
the linear span of the edge directions at $w$ has dimension $r-1$, and the multiplicities of the edges at $w$ 
have no nontrivial common factor.

\end{rem}

\medskip

Combining the previous two results and the discussion of tropical
hypersurfaces in \parref{par:tropicalization}, we obtain the following.
Recall that a {\em leaf} in a graph is a vertex of valence one.

\begin{cor}
\label{cor:unimodularisometry}
Suppose that $g(\Xhat)\geq 1$ and let $\Sigma$ be the minimal skeleton of
$\Xhat^\an$. 
\begin{enumerate}
\item If all vertices of $\Trop(X)$ are trivalent, all edges of $\Trop(X)$
  have multiplicity $1$, $\Sigma$ has no leaves,%
  \footnote{In earlier versions of this paper, the skeleton $\Sigma$ was assumed
    to be bridgeless; however, the proof only used that $\Sigma$ has no leaves.}
  and $\dim H_1(\Sigma,\RR) = \dim H_1(\Trop(X),\RR)$, then 
  $\trop : \Sigma \to \Trop(X)$ is an isometry onto its image.
\item If $X \subset \Gm^2$ is defined by a polynomial $f \in K[x,y]$ 
  whose Newton complex (see Remark~\ref{rem:newton.complex}) is a
  unimodular triangulation, 
  then $\Xhat$ has totally
  degenerate reduction and $\trop: \Sigma \to \Trop(X)$ induces an
  isometry from $\Sigma$ onto its image.
\end{enumerate}
\end{cor}

\pf Let $g' = \dim H_1(\Sigma,\RR) = \dim H_1(\Trop(X),\R)$.  Choose a compact,
connected subgraph $\Sigma'\subset\Trop(X)$ with $g' = \dim H_1(\Sigma',\R)$.
By Theorem~\ref{thm:combmult1theorem}, all $w \in \Trop(X)\cap N_G$ have
multiplicity $1$, so according to Theorem~\ref{thm:SectionThm}, there is a
(unique) finite subgraph $\td\Sigma$ of $X^\an$ mapping homeomorphically and
isometrically onto $\Sigma'$.  Since $\Sigma$ is the minimal skeleton of $X^\an$
and $\Sigma$ has no leaves, any subgraph of $X^\an$ whose first homology has
dimension at least $g'$ must contain $\Sigma$.  In particular,
$\Sigma\subset\td\Sigma$, so $\Sigma$ maps isometrically onto its image, which
proves~(1).

We now prove (2).  By \parref{par:tropicalization}, $\Trop(X)$ is trivalent with
all edges of multiplicity one.  By Baker's theorem~\cite{bakers.thm} (see
also~\cite{beelen:bakers.thm}), the genus $g(\hat X)$ is bounded by
the number of internal vertices in the Newton polytope of $f$.  Since
the Newton complex is dual to the tropicalization, each internal vertex
corresponds to a region in $\R^2\setminus\Trop(X)$; from this and the above, one
sees that $g(\hat X) \leq \dim H_1(\Trop(X),\R)$.  With $\Sigma'$ as above, by
Theorem~\ref{thm:SectionThm}, there is a (unique) finite subgraph $\td\Sigma$ of
$X^\an$ mapping homeomorphically onto $\Sigma'$ via $\trop$.  Since
$\dim H_1(\Sigma,\R) \geq \dim H_1(\td\Sigma,\R) \geq g(\hat X)$ in any case, it
follows from the genus formula~\cite[Remark~4.18]{bpr:berk_curves} that
$\dim H_1(\Sigma,\R) = g(\hat X)$ and that $\Xhat$ has totally degenerate
reduction.  In particular, the minimal skeleton $\Sigma$ has no leaves.  By~(1),
$\trop: \Sigma \to \Trop(X)$ is an isometry onto its image.  \qed

\medskip

Corollary~\ref{cor:unimodularisometry} will be important for 
applications to Tate curves in \S\ref{section:ellipticcurves}.

\begin{eg}  \label{eg:genus3-faithful}
As an example where the hypotheses of Corollary~\ref{cor:unimodularisometry} are satisfied, 
consider the genus three curve $X = V(f)$ with
\[
f= t^4(x^4 + y^4 + z^4) + t^2(x^3 y + x^3 z + xy^3 + xz^3 + y^3 z + y z^3) + t(x^2 y^2 + x^2 z^2 + y^2 z^2) + x^2 yz + xy^2z + xyz^2.
\]
See Figure~\ref{fig:genus3-faithful}.

\genericfig[ht]{genus3-faithful}{
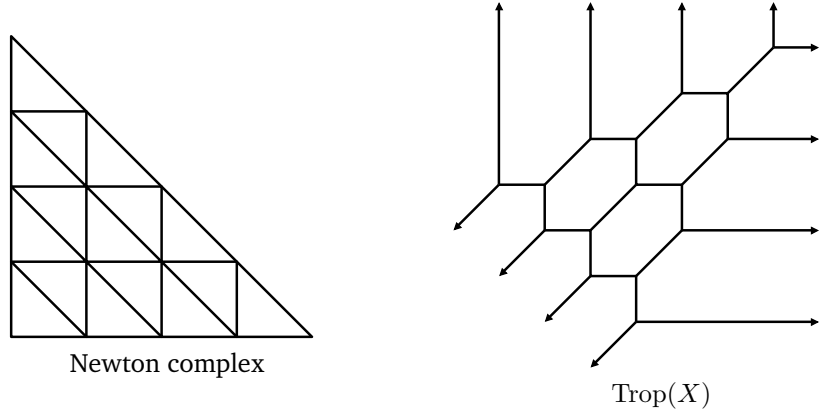
  The Newton complex and tropicalization of the curve $X$ defined by the
  polynomial $f$ from Example~\ref{eg:genus3-faithful}.  The tropicalization
  faithfully represents the minimal skeleton of $X^\berk$.}

\end{eg}

\begin{eg} \label{eg:faithful.countereg}
  The following example shows that it is possible, even under the hypotheses of Theorem~\ref{thm:SectionThm}, 
  for the tropicalization map to fail to be faithful.
  Let $K$ be the completion of $C\pu t$,
  let $\hat X \subset \bP^2$ be defined by the equation
  \[ (y-1)^2 = (x-1)^2(y+1) + t\cdot xy \]
  over $\C\pu t$, and let $X = \hat X\cap\G_m^2$.
  The above equation degenerates to a nodal rational curve when $t=0$, the
  node being $[1:1:0]$, so it defines a (not strongly) semistable
  algebraic integral model $\cX$ of $\hat X$
  --- in fact, $\cX$ is a
  minimal stable model for $\hat X$, and the associated semistable vertex
  set only contains one point, so it is minimal as well 
  (see~\cite[Corollary~4.23]{bpr:berk_curves}).
  Therefore $\hat X$ is an elliptic curve with bad reduction
  and $\inn_0(X)\cong\bar\cX\cap\G_{m,k}^2$ is reduced and irreducible.  
  The set of punctures
  \[ D \coloneq \hat X(K)\setminus X(K) = 
  \{ [0:0:1], [0:1:0], [1:0:0], [0:3:1], [2:0:1] \} \]
  reduce to distinct smooth points of $\bar\cX$, so if $\Sigma$ is the minimal
  skeleton of $\hat X^\berk$ and $\tau_\Sigma:\hat X^\berk\to\Sigma$ is the
  retraction, then $\tau_\Sigma(x)$ reduces to the generic point of
  $\bar\cX$ for all $x\in D$.  Therefore the minimal skeleton $\Gamma$ of
  $X$ and the tropicalization of $X$ are as shown in Figure~\ref{fig:faithful_countereg}.
  We see that 
  $\Trop(X)$ is contractible and everywhere multiplicity one but 
  image of the section does not contain the loop in $X^{\an}$ (the
  loop is contracted to the origin).  In particular, $\trop$ is not faithful 
  despite the fact that all points in $\Trop(X)$ have multiplicity one.
\end{eg}

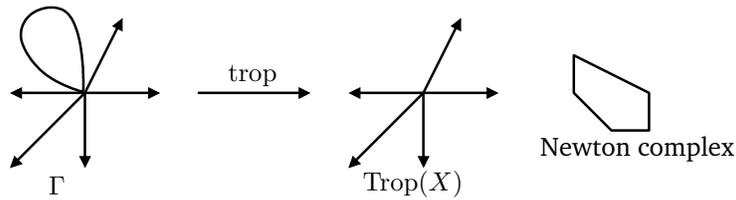
\genericfig[ht]{faithful_countereg}{The skeleton, tropicalization, and
  Newton complex of the curve $X$ from Example~\ref{eg:faithful.countereg}.  One
  sees from the Newton complex that the initial degenerations are all
  multiplicity one away from $0$, and $\inn_0(X)$ is a rational nodal
  curve.  However the tropicalization crushes the loop in $X^\berk$ to the
  origin.}

\section{Elliptic curves}
\label{section:ellipticcurves}

Let $\hat E/K$ be an elliptic curve.
If $\hat E$ has good reduction then 
the minimal skeleton $\Sigma$ of $\hat E^{\an}$ is a point, while if 
$\hat E$ has multiplicative reduction then the minimal skeleton $\Sigma$
of $\hat E^{\an}$ 
is homeomorphic to a circle of length 
$-\val(j_{\hat E}) = \val(q_{\hat E})$, where
$\hat E^\an \cong \G_m^\an / q_{\hat E}^{\ZZ}$ is the Tate uniformization
of $\hat E$ 
(see \cite[Remark 4.24]{bpr:berk_curves}).
In this section, we use our results on nonarchimedean analytic curves and their tropicalizations to prove some new results
(and reinterpret some old results) about tropicalizations of elliptic curves.

\paragraphtoc[Faithful tropicalization of elliptic curves]
As noted in \cite{KatzMarkwigMarkwig08, KatzMarkwigMarkwig09}, a curve in $\Gm^2$ given by a Weierstrass equation 
$y^2 = x^3 + ax^2 + bx + c$ cannot have a cycle in its tropicalization, because
the Newton complex of a Weierstrass equation does not have an
interior vertex.
Thus Weierstrass equations are always `bad' from the point of view of tropical geometry.
On the other hand, the following result shows that there do always exist `good' plane embeddings of elliptic curves with
multiplicative reduction.

\begin{thm}
\label{thm:goodgenus1trop}
Let $\hat E/K$ be an elliptic curve with multiplicative reduction.  Then
there is a closed embedding of $\hat E$ in $\PP^2$, given by a projective
plane equation of the form $ax^2y + bxy^2 + cxyz = dz^3$, such that
(letting $E$ be the open affine subset of $\hat E$ mapping into the torus
$\Gm^2$) $\Trop(E)$ is a trivalent graph, every point of $\Trop(E)$ has a
smooth and irreducible initial degeneration (hence tropical multiplicity
$1$), and the minimal skeleton of $\hat E$ is faithfully represented by
the tropicalization map.
In particular, $\Trop(E)$ contains a cycle of length $-\val(j_{\hat E})$.
\end{thm}

\pf Let $q = q_{\hat E}$ be the Tate parameter, so that $\hat E^\an \cong
\G_m^\an / q^{\ZZ}$.  Choose a cube root $q^{1/3} \in K^\times$ of $q$
and let $\alpha,\beta \in\hat E(K)$ correspond under the Tate
isomorphism to the classes of $q^{1/3}$ and $(q^{1/3})^2$, respectively.
Recall that a divisor $D = \sum a_P (P)$ on an elliptic curve is principal
if and only if $\sum a_P = 0$ and $\sum a_P P = 0$ in the group law on the
curve.  In particular, there exist rational functions $f$ and $g$ on 
$\hat E$
(unique up to multiplication by a nonzero constant) such that $\Div(f) =
2(\alpha) - (\beta) - (0)$ and $\Div(g) = 2(\beta) - (\alpha) - (0)$.  Let
$\psi : \hat E \to \PP^2$ be the morphism associated to the rational map
$[f:g:1]$.  Since $1,f,g$ form a basis for $L(D)$ with $D =
(\alpha)+(\beta)+(0)$ and $D$ is very ample \cite[Corollary
IV.3.2(b)]{hartshorne:ag}, $\psi$ is a closed immersion.

Let $\Gamma$ be the minimal skeleton of $E$, i.e.,\ the
smallest closed connected subset of $\hat E^{\an}$
containing the skeleton $\Sigma$ of $\hat E^{\an}$ and the three points
$\alpha,\beta,0$, with those three points removed.  Recall that $\Sigma$
is isometric to a circle of circumference 
$\ell_\an(\Sigma) = -\val(j_{\hat E}) = \val(q)$.  
The natural map
\[
K^\times \onto \hat E(K) \into \hat E^{\an} \onto \Sigma \isomap \RR / \ell \ZZ
\]
is given by $z \mapsto [\val(z)]$; in particular, if $\tau_\Sigma : \hat E^{\an} \onto
\Sigma$ denotes the canonical retraction, we have $\tau_\Sigma(0) = [0],
\tau_\Sigma(\alpha) = [\frac{1}{3} \val(q)]$, and $\tau_\Sigma(\beta) = [\frac{2}{3}
\val(q)]$.  Thus $\Gamma$ is
a circle with an infinite ray emanating from 
each of three equally spaced points $O = \tau_\Sigma(0),A = \tau_\Sigma(\alpha), B =
\tau_\Sigma(\beta)$ along 
the circle (see Figure~\ref{fig:goodtrop} below).  
The tropicalization map $\trop : E=\hat E^{\an} \smallsetminus \{ 0,\alpha,\beta \}
\to \RR^2$ corresponding to the embedding $E \into \G_m^2$ given by
$(f,g)$ factors through the retraction onto $\Gamma$.

The map $\trop : \Gamma \to \RR^2$ can be determined (up to an additive
translation) using the 
Slope Formula~\cite[Theorem~5.15]{bpr:berk_curves}
by solving an elementary graph potential problem.
The result is as follows.  The function $\val(f) = -\log|f|$ has slope $-1$
along the ray from $O$ to $0$, slope $2$ along the ray from $A$ to
$\alpha$, and slope $-1$ along the ray from $B$ to $\beta$.  On $\Sigma$,
it has slope $1$ along the segment from $O$ to $A$, slope $-1$ along the
segment from $A$ to $B$, and slope $0$ along the segment from $B$ to $O$.
Similarly, the function $\val(g) = -\log|g|$ has slope $-1$ along the ray
from $O$ to $0$, slope $-1$ along the ray from $A$ to $\alpha$, and slope
$2$ along the ray from $B$ to $\beta$.  On $\Sigma$, it has slope $0$
along the segment from $O$ to $A$, slope $1$ along the segment from $A$ to
$B$, and slope $-1$ along the segment from $B$ to $O$.  Thus (up to a
translation on $\RR^2$) $\Trop(E)$ is a trivalent graph consisting
of a triangle with an infinite ray emanating from each of the vertices as
in Figure~\ref{fig:goodtrop}.

Since the expansion factor along every edge of $\Gamma$ is equal to $1$
by~\parref{eq:mrel.is.gcd.slopes}, it follows from
Corollary~\ref{cor:Berkovichmultiplicities} that the tropical
multiplicity of every edge of $\Trop(E)$ is $1$.  By
Theorem~\ref{thm:combmult1theorem}, the multiplicity at every vertex of
$\Trop(E)$ is $1$ as well, and in fact the initial degenerations are
smooth and irreducible since the Newton complex is unimodular (see
Figure~\ref{fig:goodtrop}).  Since the expansion factor is $1$ along every
edge of $\Sigma$ and $\trop|_{\Sigma}$ is a homeomorphism, it follows that
$\Sigma$ is faithfully represented.  The bounded edges of $\Trop(E)$ form
a triangle each of whose sides has lattice length $\val(q)/3$.

The only thing which remains to be proved is that 
$\psi(\hat E) \subset \PP^2$ is cut out by an equation of the form
indicated in the statement of the 
theorem.  This follows from the Riemann-Roch theorem: the functions
$1,fg,f^2g,fg^2$ all belong to the $3$-dimensional vector space $L(3(0))$
and hence there is a nonzero linear relation between them.  (This argument
is similar to \cite[Proposition IV.4.6]{hartshorne:ag}.)\qed

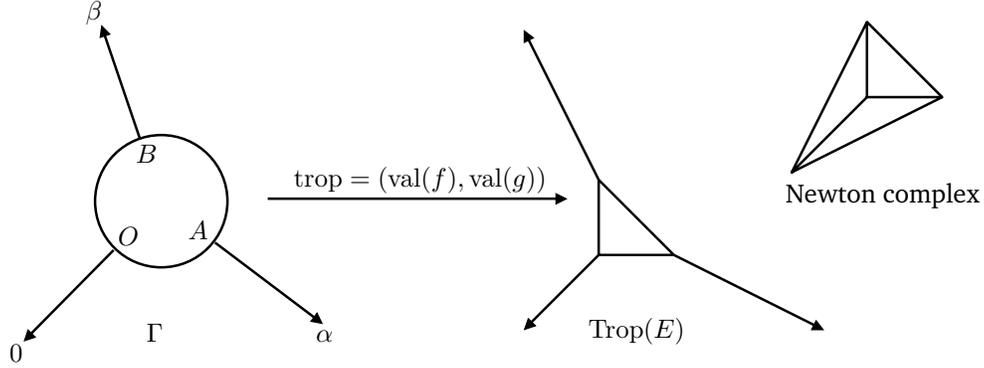
\genericfig[ht]{goodtrop}{The skeleton $\Gamma$ of $E$, the
  tropicalization of $E$, and the Newton complex of the equation
  $ax^2y + bxy^2 + cxyz = dz^3$ defining $E$, where $E$ is as in the proof
  of Theorem~\ref{thm:goodgenus1trop}.  The minimal skeleton $\Sigma$ of 
  $\hat E$ is the circle contained in $\Gamma$.  The tropicalization
  fathfully represents $\Sigma$, so 
  $\ell_{\Trop}(\Trop(\Sigma)) = \ell_\an(\Sigma)$.}

We can also use our theorems to give more conceptual proofs of many of the
results from \cite{KatzMarkwigMarkwig08, KatzMarkwigMarkwig09}.  For
example, we have the following theorem which was proved in
\cite{KatzMarkwigMarkwig09} by a brute-force computation:

\begin{thm}
\label{thm:KMM}
Let $E \subset \Gm^2$ be the intersection of an elliptic curve 
$\hat E\subset\PP^2$ with $\Gm^2$.
Assume that (i) $\Trop(E)$ contains a cycle $C$,
(ii) all edges of $\Trop(E)$ have multiplicity $1$, and (iii) $\Trop(E)$
is trivalent.  Then $\ell_{\Trop}(C) = -\val(j_{\hat{E}})$.
\end{thm}

\pf This follows immediately from Corollary~\ref{cor:unimodularisometry}(1).\qed

\begin{rem}
  Conditions (i)--(iii) from Theorem~\ref{thm:KMM} are automatically
  satisfied if the Newton complex of the defining polynomial for $E$ is a
  unimodular triangulation with a vertex lying in the interior of the
  Newton polygon, 
  as in Figure~\ref{fig:goodtrop}.  Varying the valuation of the
  coefficient corresponding the interior vertex while keeping all other
  coefficients fixed gives a natural map from an annulus in $\G_m$ to the
  $j$-line, which is finite and flat onto an annulus in the $j$-line, by
  \cite[Proposition 2.2]{bpr:berk_curves}(2).
  In particular, given a
  tropical plane curve dual to such a Newton complex and an elliptic curve
  $E$ with $j$-invariant equal to minus the length of the loop, there is
  an embedding of $E$ into a toric variety such that the tropicalization
  of the intersection with $\Gm^2$ is faithful and equal to the given
  tropical curve.  See \cite{chan_sturmfels:honeycomb} for explicit
  constructions of such embeddings for tropical curves of ``honeycomb
  normal form,'' including an algorithm for finding the honeycomb form of
  an elliptic curve, paramaterization by theta functions, the
  tropicalization of the inflection points, and relations to the group law.
\end{rem}

\begin{rem}
A different (but related) conceptual explanation for Theorem~\ref{thm:KMM} is given in \cite[Proposition 7.7]{helm_katz:monodromy}.
\end{rem}

Let us say that a closed embedding of an elliptic curve $\hat E/K$ in some
toric variety is {\em certifiably of genus $1$} if $\Trop(E)$
satisfies conditions (i)-(iii) from Theorem~\ref{thm:KMM}.  Note that the
cycle $C$ in any such embedding satisfies $\ell_{\Trop}(C) =
-\val(j_{\hat{E}})$, by Theorem~\ref{thm:KMM}.  Combining
Theorems~\ref{thm:goodgenus1trop} and \ref{thm:KMM}, we obtain:

\begin{cor}
  An elliptic curve $\hat E/K$ has multiplicative reduction if and only if it
  has a closed embedding in $\PP^2$ which is certifiably of genus $1$.
\end{cor}

\paragraphtoc[Speyer's well-spacedness condition]
In this section we explain how Speyer's well-spacedness condition
\cite{speyer:uniformizing} follows from a more general
result (possibly of independent interest) about the analytification of an elliptic curve 
$\hat E/K$.

Let $\Sigma$ be the minimal skeleton of $\hat E$.
For $P,Q \in\hat E(K)$, define $i(P,Q) \in \RR_{\geq 0} \cup \{ \infty \}$ as
follows:
\[
i(P,Q) = \left\{ 
\begin{array}{ll}
0 & {\rm \; if \; } \tau_\Sigma(P) \neq \tau_\Sigma(Q) \\
{\dist}(P \vee Q,\Sigma) &  {\rm \; if \; } \tau_\Sigma(P) = \tau_\Sigma(Q) \\
\end{array}
\right.
\]
where $\tau_\Sigma : \hat E^{\an} \onto \Sigma$ is the retraction map,
$P\vee Q$ is the first point where the geodesic paths from $P$ to $\Sigma$
and $Q$ to $\Sigma$ meet, and $\dist(x,\Sigma)$ is the distance (in
the natural metric on $\HH(\hat E^\an)$) from $x\in\hat E^\an$ to its
retraction $\tau_\Sigma(x)\in\Sigma$.  By convention we set $i(P,P) = +\infty$.  Since
translation by a point $P\in\hat E(K)$ is an automorphism of $\hat E$, it
induces an isometry on $\HH(\hat E^\an)$; therefore  $i(P,Q)$ only depends on
the difference $P-Q$ in 
$\hat E(K)$, i.e., $i(P,Q)=\iota(P-Q)$ with $\iota(R)=i(R,0)$.

The following lemma shows that $\| P, Q \| := {\rm exp}(-i(P,Q))$ is an
{\em ultrametric} on $\hat E(K)$:

\begin{lem}
\label{lem:ultrametric}
\begin{enumerate}
\item For any points $P,Q,R \in \hat E(K)$ we have $i(P,Q) \geq \min \{
  i(P,R),i(Q,R) \}$, with equality if $i(P,R) \neq i(Q,R)$.
\item If $m\in\Z$ is an integer such that $|m| = 1$ in $K$ then
  $i(mP,mQ) = i(P,Q)$ for any $P,Q\in\hat E(K)$ such that
  $i(P,Q)>0$.
\end{enumerate}
\end{lem}

\pf We begin by proving~(1).  
If either $i(P,R)=0$ or $i(Q,R)=0$ then the inequality is trivial.
Moreover, by translation invariance of $i$ we may assume that $R=0$.  So
we are reduced to showing that if $\iota(P)>0$ and $\iota(Q)>0$ then
$\iota(P-Q) \geq \min \{ \iota(P),\iota(Q) \}$.  

Let $\fE$ be the semistable formal model of
$\hat E$ corresponding to the semistable vertex set $\{\tau_\Sigma(0)\}$
(see \cite[Remark 4.24]{bpr:berk_curves}); note that
$\Sigma = \Sigma(\hat E,\{\tau_\Sigma(0)\})$).  Then $\bar\fE$ is a nodal
rational curve, and the smooth locus $\bar\fE{}^{\rm sm}$ is a group scheme
isomorphic to $\G_{m,k}$. 
The subset 
$\hat E^1(K) := \{ P \in \hat E(K) \; : \; \iota(P) > 0 \}$ is the formal
fiber over the identity element of $\bar\fE{}^{\rm sm}$, hence is a
subgroup; in fact, $\hat E^1(K)$ is isomorphic
to the group $\fm = \{ z \in K \; : \; |z| < 1 \}$ with the law of
composition given by  a one-parameter formal group law $F$ over
$R$~\cite[Proposition~VII.2.1]{silverman:I}, and the 
restriction of $\iota$ to $\hat E^1(K)$ corresponds to 
the valuation on $\fm$ under this identification.  The desired inequality
follows since a group law on $\fm$ given by a power series with
coefficients in $R$ is obviously ultrametric.

In the situation of~(2), as above we are reduced to showing that
$\iota(mP) = \iota(P)$ when $\iota(P) > 0$.  This is true because $m$ is the
coefficient of the linear term of the power series for multiplication by
$m$ under $F$~\cite[Proposition~IV.2.3]{silverman:I}, and all other terms
have larger valuation.\qed

Since $\iota(P) = \iota(-P)$, an equivalent formulation of
Lemma~\ref{lem:ultrametric}(1)  is that 
for any $P,Q \in\hat E(K)$ we have 
$\iota(P+Q) \geq \min \{ \iota(P),\iota(Q) \}$,
with strict inequality if $\iota(P) \neq \iota(Q)$.

\medskip

If $f$ is a nonconstant rational function on $\hat E$, define $N_f$ to be the
set of all $x \in\HH(\hat E^{\an})$ such that $\log|f|$ is non-constant in every
open neighborhood of $x$.  Equivalently, for $x \in\HH(\hat E^{\an})$ let $T_x(f)$
be the (finite) set of tangent directions at $x$ along which the
derivative of $\log|f|$ is nonzero.  Then $N_f$ is the set of all $x \in
\hat E^{\an}$ such that $T_x(f) \neq \emptyset$.
By~\cite[Theorem~5.15(1,2)]{bpr:berk_curves},
$N_f$ is a union of finitely many edges of the
minimal skeleton $\Gamma_f$ of the curve obtained from $\hat E^{\an}$ by
removing all zeros and poles of $f$.

\begin{thm}
\label{thm:BerkovichWellSpacedness}
Suppose that $K$ has residue characteristic zero.  Let $f$ be a nonconstant
rational function on $\hat E$ and assume that there exists $x \in N_f$ such
that ${\dist}(x,\Sigma) < {\dist}(y,\Sigma)$ for all $y \in N_f$
with $y \neq x$.  Assume also that $\Sigma\cap N_f = \emptyset$.
Then $|T_x(f)| \geq 3$.
\end{thm}

In other words, either the minimum distance from $N_f$ to the skeleton is
achieved at two distinct points, or else the minimum is achieved at a
unique point at which $\log|f|$ has nonzero slope in at least three
different tangent directions.  

\pf By the 
Slope Formula~\cite[Theorem~5.15]{bpr:berk_curves},
the sum of the outgoing slopes of $\log|f|$ at $x$ is zero, so $|T_x(f)|
\geq 2$.  Assume for the sake of contradiction that $|T_x(f)|=2$ and write
$T_x(f) = \{ v, v' \}$.  Our hypotheses imply that $x \not\in \Sigma$ and that
$\Sigma$ lies in a single connected component of $\hat{E}^{\an} \smallsetminus \{ x \}$.

Let $B(x,v)$ (resp. $B(x,v')$) be the open set
consisting of all $z \in\hat E^{\an}$ lying in the tangent direction $v$
(resp. $v'$), so that $B(x,v)$ and $B(x,v')$ are connected components of 
$\hat{E}^{\an} \smallsetminus \{ x \}$ which are disjoint from $\Sigma$.
Let $D_v$ be the restriction of $\Div(f)$ to $B(x,v)$ and let $D_{v'}$ be
the restriction of $\Div(f)$ to $B(x,v')$.  By the 
Slope Formula, 
we have $m = \deg(D_v) = -\deg(D_{v'})$ for some nonzero integer
$m$.  Without loss of generality, we may assume that $m > 0$.  Let $\delta
= {\dist}(x,\Sigma) > 0$.

We claim that $\Div(f)$ can be written as
\begin{equation}
\label{eq:divf}
\Div(f) = m\left( (P) - (Q) \right) + \sum_j \left( (A_j) - (B_j) \right)
\end{equation}
with $i(P,Q) = \delta$ and $i(A_j,B_j) > \delta$ for all $j$.  
Because we have assumed that $K$ has residue characteristic zero,
we have $|m|=1$, and therefore $i(m(P),m(Q)) = i(P,Q) = \delta$
by~Lemma~\ref{lem:ultrametric}(2).
Since $mP - mQ = \sum (B_j - A_j)$ in the group law on $\hat E(K)$, we
obtain a contradiction to Lemma~\ref{lem:ultrametric}(1).

To prove the claim, we use a trick due to D. Speyer.
Suppose $D_v = (P_1) + \cdots + (P_r) - (Q_1) - \cdots - (Q_s)$ and
$D_{v'} = (P'_1) + \cdots + (P'_{r'}) - (Q'_1) - \cdots - (Q'_{s'})$ with
$r - s = m$ and $s' - r' = m$.
Then 
\[
\begin{aligned}
D_v + D_{v'} &= m\left( (P_1) - (Q'_1) \right) + \sum_{j=1}^s \left( (P_j) - (Q_j) \right)
+ \sum_{j=s+1}^r \left( (P_j) - (P_1) \right) \\
&+ \sum_{j=1}^{r'} \left( (P'_j) - (Q'_j) \right)
+ \sum_{j=r'+1}^{s'} \left( (Q'_1) - (Q'_j) \right).
\end{aligned}
\]
Note that $i(P_1,Q_1') = \delta$ but that
$i(P_j,Q_j) > \delta$ for all $j=1,\ldots,s$, 
$i(P_j,P_1) > \delta$ for all $j=s+1,\ldots,r$,
$i(P_j',Q_j') > \delta$ for all $j=1,\ldots,r'$, and
$i(Q_1',Q_j') > \delta$ for all $j=r'+1,\ldots,s'$.

Let $C_1,\ldots,C_t$ be the connected components of $N_f$, labeled so that
$C_1$ is the component containing $x$.  By the 
Slope Formula~\cite[Theorem~5.15]{bpr:berk_curves}, 
for each $j$ the restriction $D_j$ of $\Div(f)$ to
$C_j$ is a nonzero divisor of degree zero and $\Div(f) = \sum_j D_j$.
Moreover, if $A,B \in \hat E(K) \cap C_j$ then the unique geodesic paths from
$A$ to $\Sigma$ and $B$ to $\Sigma$ must pass through the unique 
point $x_j$ of $C_j$ closest to $\Sigma$,
so for $j \geq 2$ we have $i(A,B) > \delta$.  The claim now follows since
$D_1 = D_v + D_{v'}$ can be written as above, and by what we have just
said we can write each $D_j$ for $j \geq 2$ as a sum of divisors of the
form $(A)-(B)$ with $i(A,B) > \delta$.\qed

In particular, we obtain the necessity of Speyer's well-spacedness condition for a genus $1$ tropical curve to lift:

\begin{cor}[Speyer]
\label{cor:SpeyerWellSpacedness}
Suppose that $K$ has residue characteristic zero.  Let $E$ be a dense open
subset of an elliptic curve $\hat{E}$ over $K$ with multiplicative reduction and let $\psi : E \into \T$
be a closed embedding of $E$ in a torus $\T$.  Assume that (i) every
vertex of $\Trop(E)$ is trivalent, (ii) every edge of $\Trop(E)$ has
multiplicity one, and (iii) $\Trop(E)$ contains a cycle $\Sigma'$ which is
contained in a hyperplane $H$.  If $W_H$ denotes the closure in $N_{\RR}$
of the set of points of $\Trop(E)$ not lying in $H$, then there is no
single point of $W_H$ which is closest to $\Sigma'$.
\end{cor} 

In other words, `the minimum distance from points of $\Trop(E)$ not lying in $H$ to the cycle must be achieved twice'.

\pf We may assume that $\psi : E \to \T \cong \Gm^n$ is given by
$(f_1,\ldots,f_n)$ with $\log|f_n|$ equal to a constant $c$ on $\Sigma$ 
and that $H$ is the hyperplane $x_n = c$.
Let $\Gamma$ be the minimal skeleton of $E$.  
By Corollary~\ref{cor:unimodularisometry}(1), we see that $\trop : \Gamma \to
\Trop(E)$ is an isometry.  Since $N_{f_n} \subset \Gamma$, the result now
follows from Theorem~\ref{thm:BerkovichWellSpacedness}.\qed

\begin{rem}
E. Katz \cite{katz:lifting_obstructions} and 
T. Nishinou \cite{nishinou:lifting_obstructions}
have recently obtained other kinds of generalizations of Speyer's well-spacedness condition.  Their generalized 
conditions apply to curves of higher genus.
\end{rem}

\section{A generalization of the Sturmfels-Tevelev multiplicity formula}
\label{section:drawingtropcurves}

As an illustration of the tools developed in this paper, we conclude with a generalization\footnote{A different proof is given in an appendix to \cite{osserman_payne:lifting}.} of 
the Sturmfels-Tevelev multiplicity formula
\cite[Theorem~1.1]{sturmfels_tevelev:elimination}
to the non-constant coefficient case
(and also to non-smooth points).

Let $X \subset \T$ be a closed subvariety, i.e.,\ a reduced and irreducible
closed subscheme.  Let $\alpha: \T \rightarrow \T'$ be a 
homomorphism of tori that induces a generically finite map of degree
$\delta$ from $X$ to $X'$, where $X'$ is the closure of $\alpha(X)$.
Then, set theoretically, $\Trop(X')$ is the image of $\Trop(X)$ under the
induced linear map $A: N_\R \rightarrow N'_\R$
\cite[Proposition~3]{tevelev:compactifications}.  The fundamental problem
of tropical elimination theory is to determine the multiplicities on the
maximal faces of $\Trop(X')$ from those on the maximal faces of
$\Trop(X)$.  

\begin{thm}
\label{thm:generalizedST}
Let $\alpha : \T \to \T'$ be a homomorphism of algebraic tori over $K$ and
let $X$ be a closed subvariety of $\T$.  Let $X'$ be the
schematic image of $X$ in $\T'$, let $f:X\to X'$ be the restriction
of $\alpha$ to $X$, and let $F = \trop(f): \Trop(X)\to\Trop(X')$ be the 
restriction of the linear map $A:N_\R\to N'_\R$ induced by $\alpha$.
Suppose that $f$ is generically finite of degree $\delta$.  
Then for any point $w' \in \Trop(X')\cap N'_G$ such that 
$|F^{-1}(w')|<\infty$, we have
\begin{equation}
\label{eq:ST1}
m_{\Trop}(w') = \frac{1}{\delta} \sum_{w \in F^{-1}(w')} 
\sum_{\bar\fC\subset\bar\fX\s w} \mult_{\bar\fX\s w}(\bar\fC)\,
[\bar\fC : \im(\bar\fC)],
\end{equation}
where the second sum runs over all irreducible components $\bar\fC$ of
$\bar\fX\s w$ and where $\im(\bar\fC)$ is the image of $\bar\fC$ in 
$(\bar\fX\p)\s{w'}$.
\end{thm}

In order to use the projection formula~\parref{prop:projection.formula}, we
will need the following lemma.

\begin{lem} \label{lem:degree.analytification.2}
  Let $X,X'$ be integral finite-type $K$-schemes and let
  $f:X\to X'$ be a generically finite dominant morphism of degree $\delta$.
  Let $\sU'\subset (X')\s\an$ be an analytic domain and let 
  $\sU = (f^\an)\inv(\sU')$.  If $f^\an|_\sU:\sU\to\sU'$ is finite then it 
  has pure degree $\delta$~\parref{par:pure.degree}.  
\end{lem}

\pf By~\cite[Exercise II.3.7]{hartshorne:ag}, there is a dense open subscheme
$U'\subset X'$ such that $U\coloneq f\inv(U')\to U'$ is
finite.  Shrinking $U'$ if necessary, we assume that $U'$ is smooth.
By Proposition~\ref{prop:degree.analytification} the morphism
$U^\an\to(U')\s\an$ is pure of degree $\delta$.  Let
$\sU,\sU'$ be as in the statement of the Lemma.  By
Proposition~\ref{prop:degree.cover}(2), we may assume that 
$\sU = \sM(\cA)$ and $\sU' = \sM(\cA')$ are affinoid.
By~\cite[Lemma~A.1.2(2)]{conrad:irredcomps}, $\sU$ and $\sU'$ are
equidimensional of the same dimension as $X$ and $X'$.  Therefore
$\sU'\cap(X'\setminus U')\s\an$ is nowhere dense in $\sU'$.
If $\sV' = \sM(\cB')$ is any connected affinoid subdomain of 
$\sU'\cap(U')\s\an$ then $\cB'$ is a domain because $\sV'$ is
smooth, and $\sV\coloneq (f^\an)\inv(\sV')\to\sV'$ has (pure) degree
$\delta$ because $\sV'\subset(U')^\an$.  Since 
$\sU'\cap(X'\setminus U')\s\an$ is nowhere dense in $\sU'$, we can choose
$\sV'$ such that $\Spec(\cB')\to\Spec(\cA')$ takes the generic point of 
$\Spec(\cB')$ to any given generic point of $\Spec(\cA')$.  Hence
$\sU\to\sU'$ has pure degree $\delta$.\qed

\pf[of Theorem~\ref{thm:generalizedST}]
Let $w'\in\Trop(X')\cap N'_G$ be a point with finite preimage under $F$.
Let 
\[ \sX^{w'} = \trop\inv(F\inv(w'))\cap X^\an 
= \Djunion_{w\in F\inv(w')} \sX^w. \]
This is an affinoid domain in $X^\an$ because it is a closed subspace of
the affinoid 
$\trop\inv(F\inv(w')) = \Djunion_{w\in F\inv(w')}\sU^w$.
We claim that $\sX^{w'}\to(\sX')\s{w'}$ is a finite morphism.  It suffices
to show that the composite $\sX^{w'}\to\sU^{w'}$ is a finite morphism,
where $\sU^{w'} = \trop\inv(w')\subset(\T')\s\an$.  This follows 
exactly as in the proof of Theorem~\ref{thm:finiteflatfibration}: since
$F\inv(w')$ is bounded, there is an affinoid domain of $\T^\an$
contained in $(\alpha^\an)\inv(\sU^{w'})$ and containing $\sX^{w'}$ in its
relative interior.  This means that the morphism
$\sX^{w'}\to(\sX')^{w'}$ is proper, thus finite because both spaces are
affinoid. 
Hence by Lemma~\ref{lem:degree.analytification.2} the
morphism $\sX^{w'}\to(\sX')\s{w'}$ has pure degree $\delta$.  Let
$\fX^{w'}\coloneq\Djunion_{w\in F\inv(w')}\fX^w$.  The generic fiber of
$\fX^{w'}$ is$\sX^{w'}$, and the natural morphism
$\fX^{w'}\to(\fX')\s{w'}$ is finite by
Proposition~\ref{prop:finiteness.flatness}(2) and takes generic points to
generic points by Proposition~\ref{prop:equidim.reduction}.
By the projection formula~\parref{prop:projection.formula} the induced
morphism $\bar\fX\s{w'}\to(\bar\fX\p)\s{w'}$ has pure degree
$\delta$, so summing~\eqref{eq:degree.alternate} over all irreducible
components $\bar\fC\p$ of $(\bar\fX\p)\s{w'}$ yields
\[ \delta \cdot m_{\Trop}(w') = \delta\sum_{\bar\fC\p\subset(\bar\fX\p)\s{w'}}
\mult_{(\bar\fX\p)\s{w'}}(\bar\fC\p)
= \sum_{\bar\fC\subset\bar\fX\s{w'}} \mult_{\bar\fX\s{w'}}(\bar\fC)\,
[\bar\fC:\im(\bar\fC)]. \]
Since $\fX^{w'}=\Djunion_{w\in F\inv(w')}\fX^w$, this is the
desired multiplicity formula.\qed

As a consequence of Theorem~\ref{thm:generalizedST}, we obtain:

\begin{cor}
\label{cor:STformula}
Let $\alpha : \T \to \T'$ be a homomorphism of algebraic tori over $K$ and
let $A = \trop(\alpha) : N_\R \to N'_\R$ be the natural linear map.
Let $X$ be a closed subvariety of $\T$, and suppose that $\alpha$
induces a generically finite morphism of degree $\delta$ from $X$
onto its schematic image $X'$ in $\T'$.  After subdividing, we may
assume that $A$ maps each face of $\Trop(X)$ onto a face of $\Trop(X')$.
Let $\sigma'$ be a maximal face of $\Trop(X')$.  Then
\[
m(\sigma') = \frac{1}{\delta} \sum_{A(\sigma) = \sigma'} m(\sigma) \cdot [N'_{\sigma'} : A(N_\sigma)].
\]
\noindent (Here $N_\sigma$ and $N'_{\sigma'}$ are the sublattices of $N$ and $N'$ parallel to $\sigma$ and $\sigma'$, respectively.)
\end{cor}

\pf If $w'$ is a smooth point of $\Trop(X')$ and $w$ is a smooth point of
$\Trop(X)$ with $A(w)=w'$, then $(\bar{\fX}\p)\s{w'} \cong Y' \times T'(w')$
and $\bar{\fX}\s w \cong Y \times T(w)$ with $Y,Y'$ zero-dimensional schemes
of length $m_{\Trop}(w)$ and $m_{\Trop}(w')$, respectively and
$T(w),T'(w')$ algebraic tori of dimension $\dim(X)=\dim(X')$ (cf.\
Remark~\ref{rem:split.torus}).  Moreover, 
$\alpha$ induces a finite homomorphism $T(w) \to T'(w')$ of degree
$[N'_{\sigma'} : A(N_\sigma)]$.  In
this situation, the quantity $[\bar\fC:\im(\bar\fC)]$ appearing
in~\eqref{eq:ST1} is equal to 
$[T(w):T'(w')] = [N'_{\sigma'} : A(N_\sigma)]$, and 
$m_{\Trop}(w) = \sum_{\bar\fC\subset\bar\fX\s w}\mult_{\bar\fX\s w}(\bar\fC)$,
so we are reduced to Theorem~\ref{thm:generalizedST}.\qed

\begin{rem}
The original Sturmfels-Tevelev multiplicity formula is the special case of Corollary~\ref{cor:STformula}
in which $K = k\pu T$ and $X$ is defined over $k$.
\end{rem}

\bibliographystyle{thesis}
\bibliography{BPR}
\bigskip~\bigskip

\end{document}